\documentclass[12pt,amstex]{article}

\usepackage{amsmath}
\usepackage{amsfonts}

\usepackage{amsthm}

\newtheorem{theo}{Theorem}
\newtheorem{coro}{Corollary}

\newtheorem{lemm}{Lemma}

\newcommand{\lbl}{\label}

\newcommand{\eq}[1]{$(\ref{#1})$}

\def\Bin{{\rm Bin}}
\def\N{\mathbb{N}}
\def\E{\mathbb{E}}
\def\PP{\mathbb{P}}
\def\0{0}
\def\N{\mathbb{N}}
\def\Z{\mathbb{Z}}
\def\R{\mathbb{R}}

\def\G{{\cal G}}
\def\Gn{{\cal G}(\X_n,r_n)}

\renewcommand{\E}{\mathbb E \,}

\newcommand{\tod}{\stackrel{{\cal D}}{\longrightarrow}}
\newcommand{\eqd}{\stackrel{{\cal D}}{=}}

\newcommand{\supp}{{\rm supp}}
\newcommand{\fmax}{f_{{\rm max}}}

\newcommand{\vold}{\omega_d}  

\newcommand{\eqco}{\setcounter{equation}{0}}
\newcommand{\thco}{\setcounter{theo}{0}}
\newcommand{\prco}{\setcounter{prop}{0}}
\newcommand{\laco}{\setcounter{lemm}{0}}
\newcommand{\coco}{\setcounter{coro}{0}}
\newcommand{\cjco}{\setcounter{conj}{0}}

\newcommand{\deco}{\setcounter{defn}{0}}
\newcommand{\allco}{\eqco  \thco \prco \laco \coco \cjco \deco}
\newcommand{\Po}{{\cal P}}
\newcommand{\X}{{\cal X}}
\newcommand{\bX}{{\bf X}}
\newcommand{\bx}{{\bf x}}

\newcommand{\by}{{\bf y}}

\newcommand{\polybeta}{\gamma}

\newcommand{\LL}{L}
\renewcommand{\H}{{\cal H}}
\renewcommand{\SS}{{\cal S}}
\newcommand{\TT}{{\cal T}}
\newcommand{\A}{{\cal A}}

\newcommand{\Var}{{\rm Var}}

\newcommand{\card}{{\rm card}}

\newcommand{\diam}{{\rm diam}}

\newcommand{\Y}{{\cal Y}}
\newcommand{\ZZ}{{\cal Z}}
\newcommand{\BB}{{\cal B}}
\newcommand{\U}{{\cal U}}
\newcommand{\WW}{{\cal W}}
\newcommand{\MM}{{\cal M}}

\renewcommand{\U}{{\cal U}}
\newcommand{\F}{{\cal F}}
\renewcommand{\G}{{\cal G}}

\newcommand{\NN}{{\cal N}}
\newcommand{\Cl}{\Lambda}

\newcommand{\eps}{\varepsilon}
\def\bdm{\begin{displaymath}}
\newcommand{\edm}{\end{displaymath}}
\def\benu{\begin{enumerate}}
\def\eenu{\end{enumerate}}
\def\beqn{\begin{equation}}
\def\eeqn{\end{equation}}
\def\be{\begin{equation}}
\def\ee{\end{equation}}
\def\bea{\begin{eqnarray}}
\def\eea{\end{eqnarray}}
\newcommand{\bean}{\begin{eqnarray*}}
\newcommand{\eean}{\end{eqnarray*}}
\newcommand{\bear}{\begin{eqnarray}}
\newcommand{\eear}{\end{eqnarray}}

\renewcommand{\epsilon}{\varepsilon}

\def\R{\mathbb{R}}

\def\A{{\cal A}}

\def\qed{\hfill\hbox{${\vcenter{\vbox{
    \hrule height 0.4pt\hbox{\vrule width 0.4pt height 6pt
    \kern5pt\vrule width 0.4pt}\hrule height 0.4pt}}}$}}

\def\lla{\lambda}
\def\spann{h}
\def\mmm{c}  
\def\bxdi{\gamma}  

\def\mmu{\mu}  
\def\imag{i}  

\begin{document}

\title{\bf Local central limit theorems \\ in stochastic geometry}

\author{
 Mathew D. Penrose$^1$ and Yuval Peres$^2$
\\
{\normalsize{\em
 University of Bath and
 Microsoft Research
}} }
\maketitle

 \footnotetext{ $~^1$ Department of
Mathematical Sciences, University of Bath, Bath BA1 7AY, United
Kingdom: {\texttt m.d.penrose@bath.ac.uk} }
 \footnotetext{ $~^1$ Partially supported by
the Alexander von Humboldt Foundation through
a Friedrich Wilhelm Bessel Research Award.}
\footnotetext{$~^2$  Microsoft research,  Redmond, WA USA: {\texttt peres@microsoft.com}}

\begin{abstract}
We give a general local central limit theorem
for the sum of two independent random variables, one of
which satisfies a central limit theorem while
the other satisfies a local central limit theorem
with the same  order  variance. We apply
this result to various quantities
arising in stochastic  geometry, including: size of the largest
component for percolation on a box;
 number of components, number of edges, or number
of isolated points, for random geometric graphs; covered volume
for germ-grain coverage models; number of accepted points for
finite-input random sequential adsorption; sum of nearest-neighbour
distances for a random
sample from a continuous
multidimensional
distribution.
\end{abstract}



\newpage

\tableofcontents

\section{Introduction}
\lbl{secintro}
A number of general central limit theorems (CLTs) have
been proved recently for quantities arising in stochastic  geometry
 subject to a certain local dependence.
  See \cite{Penbk,PenCLT,PenEJP,PY1,PY2} for some examples.
The present work is concerned with
{\em local}~ central limit theorems
for such quantities.
The local CLT for a binomial $(n,p)$ variable
says that
for  large $n$ with $p$ fixed, its probability mass
function
minus that of the corresponding normal variable
rounded to the nearest integer,
is uniformly $o(n^{-1/2})$.
The classical local CLT provides
 similar results  for sums of i.i.d.\ variables with an arbitrary distribution
possessing a finite second moment.
 Here we are concerned with
sums of variables with some weak dependence, in the sense that
the summands can be thought of as contributions  from spatial
regions with only local interactions between different regions.

Among the examples for which we obtain local CLTs here are the following.
In  Section \ref{secperc} we give local CLTs for
the number of clusters in percolation on a large finite lattice box,
and for the size of the largest open cluster for supercritical
percolation on a large finite box, as the box size becomes large.
In Sections \ref{secRGG} and \ref{secstogeo}
we consider continuum models, starting with random geometric
graphs \cite{Penbk} for which we demonstrate local CLTs for the
number of copies of a fixed subgraph (for example the number of
edges) both in the thermodynamic limit 
 (in which the mean degree is $\Theta(1)$)
 and in the sparse limit  (in which the mean degree vanishes).
For the thermodynamic limit we also  derive local CLTs for
the number of components of a given type (for example the number
of isolated points),  as an example of a more general local CLT
for functionals which have finite range interactions
or which are sums of functions determined by nearest neighbours
(Theorem \ref{finranthm}).  This
also yields  local CLTs for quantities associated
with a variety of other models, including germ-grain models and
random sequential adsorption in the continuum.

We derive these  local CLTs
using the following idea which has  been seen (in somewhat different
form) in \cite{DMcD}, in \cite{Be},
and no doubt
 elsewhere.
If the random variable of interest
is known to satisfy a CLT, and can be decomposed
(with high probability) as the sum
of two independent parts, one of which satisfies a local CLT with
the same order of  variance growth, then one can find a local CLT for the
original variable.
Theorem \ref{genthm2} below formalises this idea. The statement of
this result has no geometrical content and it could be of use elsewhere.

 In the geometrical context, one can often use the geometrical
structure to effect such a decomposition. Loosely speaking,
in these examples one can represent a positive proportion
of the spatial region under consideration as a union of
disjoint boxes or balls, in such a way that with high probability a
non-vanishing  proportion of the boxes are `good' in some sense,
where the contributions to the variable of interest from a good box,
given the configuration outside the box and given that it has
the `good' property, are i.i.d.
Then the classical
local CLT applies to the total contribution from good boxes,
and one can represent
the variable of interest as the sum of two independent contributions,
one of which (namely the contribution from good boxes)
 satisfies a local CLT, and then apply Theorem \ref{genthm2}.
This technique is related to a method used by Avram and Bertsimas
\cite{AB} to find lower bounds on the variance for certain
quantities in stochastic geometry, although the examples considered
here are mostly different from those considered in \cite{AB}.

In any case, our results provide extra information
on the CLT behaviour for variables for numerous geometrical and
multivariate stochastic settings, which have arisen in  a variety
of applications (see the examples in Section \ref{secstogeo}).

\section{A general local CLT}
\lbl{secgenresult}
In the sequel we let $\phi$ denote the standard
($\NN(0,1)$)
normal density function, i.e. $\phi(x)= (2\pi)^{-1/2}
\exp(-(1/2)x^2)$. Note that for $\sigma >0$,
 the probability density function of the
$\NN(0,\sigma^2)$
 distribution
is then $\sigma^{-1}\phi(x/\sigma)$, $x \in \R$.
Define the $\NN(0,0)$ distribution
to be that of a random variable that
is identically zero.

We say a random variable $X$ 
is {\em integrable} if $\E[|X|]< \infty$. We say
$X$ has a {\em lattice} distribution if there
exists  $h >0$
 such that   $(X-a)/h \in \Z$ almost surely
 for some $a \in \R$.
If $X$ is lattice,  then the largest
such $h$ is called the {\em span} of $X$, and here
denoted $\spann_X$. If $X$ is non-lattice, then we
set $\spann_X :=0$.  If $X$ is degenerate, i.e.
if $\Var[X]=0$,  then we set
$\spann_X := +\infty$.
As usual with local central limit theorems, we need to
distinguish between the lattice and non-lattice cases.
For real numbers $a \geq 0,b >0$, we shall write $a|b$
to mean that  either
 $b$ is an integer multiple of $a$ or $a =0$.
When $a=  +\infty, b< \infty$
we shall say by convention that $a|b$ does not hold.

\begin{theo}
\lbl{genthm2}
Let
 $V,V_1,V_2, V_3,\ldots$ be independent identically
distributed random variables.
  Suppose for each $n \in \N$ that
$(Y_n,S_n,Z_n)$ is a triple of integrable
 random variables on the same sample space such that
(i) $Y_n$ and $S_n$ are independent, with
$S_n \eqd \sum_{j=1}^n V_j$;
(ii) both
$ n^{-1/2} \E[|Z_n - (Y_n+ S_n)|]$ and
$ n^{1/2} P[Z_n \neq Y_n+ S_n]$ tend to zero as
$n \to \infty$;
and (iii) for some $\sigma \in [0,\infty)$,
\bea
n^{-1/2} (Z_n - \E Z_n ) \tod \NN(0,\sigma^2)
{\rm ~~~as ~~} n \to \infty.
\lbl{normlim1a}
\eea
Then $\Var [ V ] \leq \sigma^2$
and if $b, \mmm_1, \mmm_2, \mmm_3, \ldots$ are positive constants with $\spann_V | b$
and $\mmm_n \sim n^{1/2}$ as $n \to \infty$, then
\bea
\sup_{u \in \R }
\left\{ \left|
\mmm_n P[Z_n \in [u,u +b) ] -
\sigma^{-1} b \phi \left(\frac{u - \E Z_n}{ \mmm_n
 \sigma
}
\right)
\right|
\right\}
\to 0
~~~~
{\rm as} ~ n \to \infty.
\lbl{1102c2}
\eea
Also, 
\bea
n^{-1/2} (Y_n- \E Y_n) \tod \NN(0,\sigma^2 - \Var [V ] ).
\lbl{0110a}
\eea
\end{theo}
\noindent{\bf Remarks.} The main case to consider is $\mmm_n=n^{1/2}$. The more general formulation above is convenient in some applications, e.g., in the proof of Theorem \ref{Gthm}.
  Theorem \ref{genthm2} is proved in Section \ref{secpfgen}.
Our main interest is in the conclusion \eq{1102c2},
 but \eq{0110a}, which comes out
for free from the proof, is also of interest.

\section{Percolation}
\lbl{secperc}
\allco
Most of our applications of Theorem \ref{genthm2} will be in the continuum,
but we start with applications to percolation on the lattice.
We consider {\em site percolation} with parameter $p$, where
each site (element) of $\Z^d$ is open with probability $p$ and closed otherwise,
independently of all the other sites.  Given a finite set
 $B \subset \Z^d$,
the {\em open clusters in $B$} are
defined to be
the components of the (random) graph with vertex set
consisting of the open sites in $B$, and edges
between each pair of open sites in $B$ that are
 at unit Euclidean distance from each other.
Let $\Cl(B)$ denote the {\em number of open clusters} in $B$.
Listing the open clusters in $B$ as ${\bf C}_1,\ldots,{\bf C}_{\Cl(B)}$,
and denoting by $|{\bf C}_j|$ the order
(i.e., the number of vertices)
of  the cluster ${\bf C}_j$, we
denote by $\LL(B)$ the random variable $\max(|{\bf C}_1|, \ldots,
|{\bf C}_{\Cl(B)}|)$, and refer to this as the {\em size of
the largest open cluster in $B$}.
Given a growing sequence of regions $(B_n)_{n \geq 1}$ in $\Z^d$,
 we shall demonstrate local CLTs for
the random variables $\Cl(B_n)$ and $\LL(B_n)$,
 subject to some conditions on the sets $B_n$
which are satisfied, for example, if they are cubes of side $n$.
 There should not be any difficulty adapting these
results to bond percolation.

For $B \subset \Z^d$ let  $|B|$ denote the number of elements
of $B$. Let $| \partial B|$ denote the number of elements of $\Z^d \setminus B$
lying at unit Euclidean distance from some element of $B$.
We say a sequence $(B_n)_{n \geq 1}$ of nonempty finite sets in $\Z^d$
has
 {\em vanishing relative boundary}
 if 
 \bea
\lim_{n \to \infty}
|\partial B_n | /|B_n| =
 0.
\lbl{vrb}
\eea
 We write $\liminf(B_n)$ for
$\cup_{n \geq 1} \cap_{m \geq n} B_m$.

\begin{theo}
\lbl{thclus}
Suppose $d \geq 2$ and $p \in (0,1)$.
Then there
exists $\sigma >0$ such that if $(B_n)_{n \geq 1}$ is
any sequence of nonempty finite subsets in $\Z^d$ with vanishing
relative boundary and with $\liminf(B_n)= \Z^d$,
then
\bea
|B_n|^{-1/2} (\Cl(B_n) - \E \Cl(B_n) ) \tod \NN(0,\sigma^2)
\lbl{ClusCLT}
\eea
and
\bea
\sup_{j \in \Z }
\left | |B_n|^{1/2} P [\Cl(B_n) = j ] - \sigma^{-1} \phi  \left(
\frac{j- \E \Cl(B_n)}{\sigma |B_n|^{1/2}} \right) \right| \to 0.
\lbl{ClusLLT}
\eea
\end{theo}

For the size of the largest open cluster we consider a more restricted
class of sequences $(B_n)_{n \geq 1}$.
Let us say
that
$(B_n)_{n \geq 1}$
is a {\em cube-like sequence of lattice boxes} if each set $B_n$
is of the form $\prod_{j=1}^d([-a_{j,n},b_{j,n} ] \cap \Z)$, where
$a_{j,n} \in \N $ and $b_{j,n} \in \N$ for all $j,n$, and moreover
\bea
\liminf_{n \to \infty} \frac{\inf \{a_{1,n},b_{1,n},a_{2,n},b_{2,n},
\ldots,a_{d,n},b_{d,n} \}}{
\sup \{a_{1,n},b_{1,n},a_{2,n},b_{2,n},
\ldots,a_{d,n},b_{d,n} \}} >0
\lbl{cubelike}
\eea
which says, loosely speaking, that the sets $B_n$ are not too far away from
all being cubes.

Given $d \geq 2$, and $p \in (0,1)$, let $\theta_d(p)$ denote the
 percolation probability, that is, the probability
that the graph with vertices consisting of all open sites in
$\Z^d$ and edges between any two open sites that are  unit Euclidean
distance apart includes an infinite component containing the origin.
Let $p_c(d)$ denote the critical value of $p$ for site percolation
in $d$ dimensions, i.e., the infimum of all $p \in (0,1)$ such
that $\theta_d(p) >0$.
It is well known that $p_c(d) \in (0,1) $ for all $d \geq 2$.

\begin{theo}
\lbl{thlargest}
Suppose $d \geq 2$ and $p \in (p_c(d),1)$.
Then there
exists $\sigma >0$ such that if $(B_n)_{n \geq 1}$ is
any cube-like sequence of lattice boxes $\Z^d$ with
 $\liminf(B_n)= \Z^d$,
we have
\bea
|B_n|^{-1/2} (\LL(B_n) - \E \LL(B_n) ) \tod \NN(0,\sigma^2)
\lbl{LargCLT}
\eea
and
\bea
\sup_{j \in \Z }
\left | |B_n|^{1/2} P [\LL(B_n) = j ] - \sigma^{-1} \phi  \left(
\frac{j- \E \LL(B_n)}{\sigma |B_n|^{1/2}} \right) \right| \to 0.
\lbl{LargLLT}
\eea
\end{theo}
Theorems \ref{thclus} and \ref{thlargest}
are proved in Section \ref{secpfperc}.
Theorem \ref{thclus} is the simplest of our applications of Theorem
\ref{genthm2}
and we give its proof with some extra detail for instructional purposes.

\section{Random geometric graphs}
\lbl{secRGG}
\allco
For our results in this section and the next,
 on continuum stochastic geometry,
let 
 $X_1,X_2,\ldots$ be i.i.d. $d$-dimensional random  vectors
 with common  density $f$.
Assume throughout that $\fmax := \sup_{x \in \R^d} f(x) < \infty$,
and that $f$ is almost everywhere continuous.
Define the induced {\em binomial point processes}
\begin{equation} \label{bin}
\X_n:= \X_{n}(f): = \{X_{1},...,X_{n}\}, ~~~ n \in \N.
\end{equation}
In the special case where $f$  is the density of the uniform distribution
on the unit $[0,1]^d$ cube we write $f \equiv f_U$.

For locally finite $\X \subset \R^d$ and $r >0$,
let $\G(\X,r)$ denote the graph with vertex set $\X$ and
with edges connecting each pair of vertices $x,y$ in $\X$ with
$|y-x| \leq r$; here $|\cdot|$ denotes the Euclidean norm though there
should not be any difficulty extending our results to other norms.
 Sometimes $\G(\X,r)$ is called a {\em geometric graph} or
 {\em Gilbert graph}.

Let $(r_n)_{n \geq 1}$ be a sequence with $r_n \to 0$ as $n \to \infty$.
  Graphs of the type of $\Gn$ are the
subject of the monograph \cite{Penbk}.
Among the quantities  of interest
associated with $\Gn$  are the number
of edges, the number of triangles, and so on;  also
the number of isolated points, the number of isolated edges, and so on.
CLTs for such quantities are given in
Chapter 3 of \cite{Penbk} (see the notes therein for other
references) for a
large class of limiting regimes for $r_n$.
Here we give some associated local CLTs.

Let $\kappa \in \N$ and let $\Gamma$ be a
fixed connected graph with $\kappa$ vertices. 
We follow terminology in \cite{Penbk}.  With $\sim$ denoting
graph isomorphism,  let $G_n$ be
the number of $\kappa$-subsets $\Y$ of $\X_n$
such that $\G(\Y,r_n) \sim \Gamma$
(i.e., the number of induced subgraphs of $\Gn$ that are isomorphic
to $\Gamma$). Let
$G^*_n$ (denoted
$J_n$ in \cite{Penbk}) denote the number of {\em components} of $\Gn$ that
are isomorphic to $\Gamma$. To avoid
certain trivialities, assume that $\Gamma$ is {\em feasible}
in the sense of \cite{Penbk}, i.e. that $\G(\X_\kappa,r)$
is isomorphic to $\Gamma$ with strictly positive probability for
some $r >0$. When considering $G_n$,
we shall also assume that $\kappa \geq 2$.
We shall give local CLTs for $G_n$ and $G^*_n$.

We assume existence of the limit
\bea
\rho :=
\lim_{n \to \infty} (n r_n^d)
< \infty ,
\lbl{rhofin}
\eea
so that $\rho $ could be zero.
If $\rho > 0$ then we are taking the {\em thermodynamic limit}.

We also assume that
\bea
\tau_n^2 := n (n r_n^d)^{\kappa -1} \to \infty ~~~ {\rm as}
~~
n \to \infty.
\lbl{taubig}
\eea
Then (see Theorems 3.12 and 3.13 of \cite{Penbk})
 there exists a constant $\sigma = \sigma(f,\Gamma,\rho) >0$,
given explicitly in terms of $f,\Gamma$ and $\rho$ in \cite{Penbk},
 such that
\bea
\lim_{n \to \infty} \tau_n^{-2} \Var( G_n) = \sigma^2;
\lbl{varlim1}
\\
\tau_n^{-1}
(G_n - \E G_n)
 \tod N(0,\sigma^2).
\lbl{normlim1}
\eea
We prove here an associated  local central limit theorem for the case $f
\equiv f_U$.
\begin{theo}
\lbl{Gthm}
Suppose $f \equiv f_U$.
Suppose $k \geq 2 $, and suppose assumptions \eq{rhofin} and \eq{taubig} hold.
Then as $n \to \infty$,
\bea
\sup_{j \in \Z }
\left | \tau_n P [G_n = j ] - \sigma^{-1} \phi  \left(
\frac{j- \E G_n}{\sigma \tau_n} \right) \right| \to 0.
\lbl{LLT_Ga}
\eea
\end{theo}
 We prove Theorem \ref{Gthm} in Section \ref{secpfedges}.
It should be possible to obtain similar results  for $G^*_n$, but we
shall do so only for the thermodynamic limit with $\rho  >0$,
 as an example in the next section.  In the next section we shall see that
for the case with $\rho >0 $, it
 is possible to relax the assumption that $f\equiv f_U$ in Theorem \ref{Gthm};
when $\rho =0$,
a similar extension to non-uniform densities should be possible, but
we content  ourselves
here with the case $f \equiv f_U$ so as to provide one example where the
simplicity and the appeal of the approach
 do not get buried.

\section{General local CLTs in stochastic geometry}
\lbl{secstogeo}
\allco

In this section we present some general local central limit theorems
in stochastic geometry. We shall illustrate these by
some examples in the next section.

For our general local CLTs in stochastic geometry,
 we consider {\em marked}
point sets in $\R^d$. Let $\MM$ be an arbitrary measurable
space (the {\em mark space}), and let $\PP_\MM$ be
a probability distribution on $\MM$. 
Given $\bx = (x,t) \in \R^d \times \MM$ and
given $y \in \R^d$, set $ y + \bx := (y+x,t)$.
Given also $a \in \R$, set $a \bx = (ax,t)$.
We think of $t$ as a mark attached to the point $x\in \R^d$
that is unaffected by translation or scalar multiplaction.
Given $\X^* \subset \R^d \times \MM, y \in \R^d$,
 and
$a \in (0,\infty)$, let $y + a\X^* := \{y + a\bx: \bx \in \X^* \}.$
Let $\0$ denote the origin of $\R^d$.
 For $x \in \R^d$,
and $r>0$, let $B(x;r)$ denote the
Euclidean ball $\{ y \in \R^d: |y-x| \leq r\}$,
and set $B^*(x;r) := B(x;r) \times \MM$.
Set $B(r) := B(\0;r)$ and $B^*(r):=B^*(0;r)$.
Given nonempty $\X^* \subset \R^d \times \MM$ and $ \Y^* \subset \R^d \times \MM$,
write
 $$
D(\X^*,\Y^*):= \inf\{|x-y|:(x,t) \in \X^*, (y,u) \in \Y^* ~~{\rm for~ some~}
t,u \in \MM\}.
$$
Let $\vold$ denote the volume of the
$d$-dimensional unit ball $B(1)$.

Suppose $H(\X^*)$ is a  measurable $\R$-valued function
defined for all
finite  $\X^* \subset \R^d \times \MM$.
 Suppose $H$ is translation invariant,
i.e. $H(y+\X^*)= H(\X^*)$ for all $y \in \R^d$ and all
$\X^*$.

Throughout this section we consider the thermodynamic limit;
let $r_n, n \geq 1  $ be a sequence of constants
such that \eq{rhofin} holds with $\rho >0$.
Define
\bea
H_n(\X^*) := H( r_n^{-1} \X^* ).
\lbl{0110b}
\eea
Let the point process $\X_n := \{X_1, \ldots,X_n\} $ in $\R^d$
 be as given in \eq{bin},
with $f$ as in Section \ref{secRGG} (so $\fmax < \infty$ and
$f$ is Lebesgue-almost everywhere continuous). Define 
the corresponding marked point processs (i.e., point
process in $\R^d \times \MM$) by
\bean
\X_n^* := \{(X_1,T_1),\ldots,(X_n,T_n)\} ,
\eean
where $(T_1,T_2,T_3, \ldots)$ is a sequence of  independent $\MM$-valued
random variables with distribution $\PP_\MM$, independent
of everything else.
We are interested in local CLTs for $H_n(\X^*_n)$, for
 general functions $H$.
 We give two distinct
types of condition on $H$, either of which is sufficient
to  obtain a local CLT.

We shall say that $H$ has {\em finite range  interactions}
if there exists a constant   $\tau  \in (0,\infty)$ such that
\bea
H(\X^* \cup \Y^*) = H(\X^*) + H(\Y^*)  ~~~ {\rm whenever}
~~
D(\X^*,\Y^*) > \tau.
\lbl{finraneq}
\eea

In many examples it is natural to write $H(\X^*)$ as a sum.
Suppose $\xi(\bx; \X^*)$ is a  measurable $\R$-valued function
defined for all pairs $(\bx,\X^*)$, where $\X^* \subset \R^d\times \MM$
 is finite  and $\bx$
is an element of $\X^*$.
 Suppose $\xi$ is translation invariant,
i.e. $\xi(y+\bx;y+\X^*)= \xi(\bx;\X^*)$ for all $y \in \R^d$ and all
$\bx,\X^*$.  Then $\xi$   induces a
translation-invariant functional $H^{(\xi)}$ defined on finite
point sets $\X^* \subset \R^d \times \MM$ by
  \bea H^{(\xi)}(\X^*) :=
\sum_{\bx \in \X^*} \xi(\bx; \X^*).
\label{induceh}
\eea
Given $r \in (0,\infty)$ we say $\xi$
has  range  $r$ if  $\xi((x,t); \X^*) = \xi ((x,t); \X^* \cap B_r^*(x))$ 
for all  finite $\X^* \subset \R^d \times \MM$ and all $ (x,t) \in \X^*$.
It is easy to see that if $\xi$ has range $r$ for some (finite) $r$ then
 $H^{(\xi)}$ has  finite range interactions,
although not all  $H$ with finite range interactions
arise in this way.

Let $\kappa \in \N$.
Given any set $\X^* \subset \R^d \times \MM$
with more than $\kappa$ elements, and given $\bx =(x,t) \in \X^*$, set
$R_\kappa(\bx;\X^*)$ to be the $\kappa$-nearest neighbour distance
from $x$ to $\X^*$, i.e. the smallest $r \geq 0$ such that
$\X^* \cap B^*(x;r)$ has at least $\kappa$ elements other than
$\bx$ itself. If $\X^*$ has $\kappa$ or fewer elements, set
$R_\kappa(\bx;\X^*):= \infty$.

We say that $\xi$ {\em depends only on the $\kappa$ nearest neighbours} if
for all $\bx$ and $\X^*$, writing $\bx = (x,t)$ we have
$$
\xi(\bx;\X^*) = \xi(\bx ; \X^* \cap B^*(x; R_\kappa(\bx;\X)) ).
$$
We give local CLTs
for $H$ under two alternative sets of conditions: either (i)
when $H$ has finite range interactions, or (ii) when $H$ is induced,
according to the definition \eq{induceh},
 by a functional $\xi(\bx;\X^*)$ which depends only on the $\kappa$
nearest neighbours, for some fixed $\kappa$.

Given $K>0$ and $n \in \N$, define point processes
$\U_{n,K},$ and 
$\ZZ_n$ 
in $\R^d$, and 
point processes
$ \U^*_{n,K}$,  
and $\ZZ_n^*$ 
in $\R^d \times \MM$, as follows.
Let $\U_{n,K}$ denote the  point process consisting of
$n$ independent uniform random points $U_{1,K},\ldots, U_{n,K}$
in $B(K)$, and let
 $\ZZ_n$ be the point process consisting of $n$ independent
points $Z_1,\ldots,Z_n$ in $\R^d$, each with a $d$-dimensional standard normal
distribution (any other positive continuous density on
$\R^d$ would do just as well). 
The corresponding marked point processs
are defined by
\bean
\U_{n,K}^* := \{(U_{1,K},T_1),\ldots,(U_{n,K},T_n)\} ;
\\
\ZZ_n^* := \{(Z_1,T_1),\ldots,(Z_n,T_n)\}.
\eean
Define the limiting span 
\bea
 \spann(H) := \liminf_{n \to \infty} \spann_{H(\ZZ^*_n)}.
 \lbl{limspandef}
\eea

\begin{theo}
\lbl{finranthm}
 Suppose that either (i)  $H$ has finite range interactions
and $h_{H(\ZZ_n^*)} < \infty$ for some $n \in \N$,
or (ii) for some $\kappa \in \N$, 
$H$ is induced by a functional $\xi(\bx;\X^*)$
which depends only on the $\kappa$ nearest neighbours,
and $h_{H(\ZZ^*_n)} < \infty$ for some $n \in \N$
with $n > \kappa$.  Suppose also that $H_n(\X_n^*)$ and
$H(\U^*_{n,K})$ are integrable for all $n \in \N$ and $K > 0$.
Finally suppose that 
\bea
n^{-1/2} ( H_{n}(\X^*_n) - \E H_{n}(\X^*_n) ) \tod \NN (0,\sigma^2) 
~~~~ {\rm as} ~ n \to \infty.
\lbl{Hclteq}
\eea
  Then
 $\sigma >0$ and
 $\spann(H) < \infty$ and
for any $b \in (0,\infty)$, with $\spann(H)| b$,
\bea
\sup_{u \in \R }
\left\{ \left|
 n^{1/2} P[H_n (\X^*_n) \in [u,u +b) ] -
\sigma^{-1} b \phi \left(\frac{u - \E H_n(\X^*_n)}{ n^{1/2} \sigma}
\right)
\right|
\right\}
\to 0
~~~~
{\rm as} ~ n \to \infty.
\nonumber \\
\lbl{1205a}
\eea
\end{theo}
We prove Theorem \ref{finranthm} in Section \ref{secpfSG}.
Analogues to this result and  to Theorem \ref{Gthm} should
also hold if one Poissonizes the number of points in the
sample, but we do not give details.

The corresponding result for unmarked point sets in
$\R^d$ goes as follows; we adapt our terminology to
this case in an obvious manner.
\begin{coro}
\lbl{finrancoro}
Suppose $H(\X)$ is $\R$-valued and defined for all
finite $\X \subset \R^d$. Suppose $H$ is translation
invariant, and set $H_n(\X):= H(r_n^{-1} \X)$.
 Suppose that either (i)  $H$ has finite range interactions
and $h_{H(\ZZ_n)} < \infty$ for some $n \in \N$,
or (ii) for some $\kappa \in \N$, $H$ is induced by a functional 
$\xi(x;\X)$
which depends only on the $\kappa$ nearest neighbours,
and $h_{H(\ZZ_n)}< \infty$ for some $n \in \N$ with $n > \kappa$. 
Suppose also that $H_n(\X_n)$ and
$H(\U_{n,K})$ are integrable
for all $n \in \N$ and $K > 0$. 
Finally suppose
\bea
n^{-1/2} ( H_{n}(\X_n) - \E H_{n}(\X_n) ) \tod \NN (0,\sigma^2) 
~~~ {\rm as} ~ n \to \infty.
\lbl{0623a}
\eea
Then $\sigma >0$ and $\spann(H) < \infty$ and 
for any $b \in (0,\infty)$, with $\spann(H) | b$,
\bea
\sup_{u \in \R }
\left\{ \left|
 n^{1/2} P[H_n (\X_n) \in [u,u +b) ] -
\sigma^{-1} b \phi \left(\frac{u - \E H_n(\X_n)}{ n^{1/2} \sigma}
\right)
\right|
\right\}
\to 0
~~~~
{\rm as} ~ n \to \infty.
\nonumber \\
\lbl{0623b}
\eea
\end{coro}
Corollary \ref{finrancoro} is easily obtained
from Theorem \ref{finranthm}
by taking $\MM$ to have just a single element, denoted $t_0$ say, 
and identifying each element $(x,t_0)\in \R^d \times \MM$
with the corresponding element $x$ of $\R^d$.

To apply Theorem \ref{finranthm} in examples, we need
to check condition \eq{Hclteq}.
For some examples this is best done directly.
However, if we strengthen the other 
hypotheses of Theorem 5.1, we can obtain
\eq{Hclteq} from known results and so do not
need to include it as an extra hypothesis.
 The next three theorems illustrate this. 
As well \eq{Hclteq}, these results 
 give us the associated variance
convergence result
\bea
\lim_{n \to \infty} n^{-1} {\rm Var} [ H_n(\X_n^*) ] = \sigma^2. 
\label{varconv}
\eea
In the next three theorems, we impose
some extra assumptions besides those of
 Theorem \ref{finranthm}. 
Writing $\supp(f)$ for the support of $f$,
we shall assume that $\supp(f)$ is compact, and that
also $r_n$ satisfy 
\bea
\lbl{strongrn}
|r_n^{-d} -n| = O(n^{1/2}),
\eea
which implies \eq{rhofin} with $\rho =1$.
We also assume certain polynomial growth bounds; see
\eq{polyxi}, \eq{polybd} and \eq{polynbr} below.

First consider the case
 where $H = H^{(\xi)}$ is induced by a functional 
$\xi(\bx;\X^*)$ with finite range $r >0$.
For any set $A$, let  $\card(A)$ denotes the number of elements of
$A$.

\begin{theo}
\lbl{fin0theo}
Suppose $H = H^{(\xi)} $ is induced by a translation
invariant functional $\xi(\bx;\X^*)$ having  finite range $r$ and
and satisfying  for some $\polybeta >0$ the polynomial growth bound
 \bea
 |\xi((x,t);\X^*) | \leq \polybeta  (\card ( \X^* \cap B^*(x;r)))^\polybeta ~~~
  \forall ~{\rm finite}~
 \X^* \subset \R^d \times \MM,  ~ \forall ~ (x,t) \in \X^*.   
 \nonumber \\
 \lbl{polyxi}
 \eea
Suppose $h_{H(\ZZ^*_n)} < \infty$ for some $n \in \N$,
 and suppose $\supp(f)$ is
 compact.
Finally, suppose that
\eq{strongrn} holds.
Then there exists $\sigma \in (0,\infty)$
 such that  \eq{Hclteq}  and \eq{varconv} hold, and   
$\spann(H) < \infty$ and
\eq{1205a} holds for all $b$ with $\spann(H)|b$.
\end{theo}

 Now we turn to the general case of Condition (i)
in Theorem \ref{finranthm}, where $H$ has finite range interactions
but is not induced by a finite 
range $\xi$. For this case we shall borrow some concepts from
continuum percolation.
  For $\lla >0$,
let $\H_\lla$ denote a homogeneous Poisson point process
in $\R^d$ with intensity $\lla$.
Let $\H_{\lla}^*$ denote the same Poisson point process
with each point given an independent $\MM$-valued mark
with the distribution $\PP_\MM$.

  Let $\lla_c$ be the critical
value for percolation in $d$ dimensions, that is, the
supremum of the set of all $\lla >0$ such that
the component of the geometric (Gilbert) graph
$G(\H_\lla \cup\{\0\},1)$
containing the origin is almost surely finite. It is known (see
e.g. \cite{Penbk})
that $0 < \lla_c < \infty$ when $d \geq 2$ and $\lla_c =  \infty$
when $d =1$.

For nonempty $\X \subset \R^d$, write $\diam(\X)$ for
 $\sup \{ |x-y|: x,y \in \X\}$.
For  $\X^* \subset \R^d \times \MM$, write $\diam(\X^*)$
for $\diam (\pi(\X^*))$,
where $\pi$ denotes the canonical projection from
$\R^d \times \MM $ onto $\R^d$. 

\begin{theo}
\lbl{fintheo}
Suppose $H(\X^*)$ is a measurable $\R$-valued function defined for
all finite $\X^* \subset \R^d \times \MM$, and is translation invariant.
Suppose $\supp(f)$  is compact.  Suppose for some $\tau >0$
that the finite range interaction condition 
  \eq{finraneq} holds, and suppose $f$ and $  \tau$ satisfy
the subcriticality condition
\bea
 \tau^d \fmax < \lla_c,
\lbl{subcrit}
\eea
Assume $(r_n)_{n \geq 1}$ satisfies \eq{strongrn},
and suppose also that
 $ h_{H(\ZZ^*_n)}  <\infty$ for some $n \in \N$, 
and that there exists a constant $\polybeta >0$
such that for all finite non-empty $\X^* \subset \R^d$ we have
\bea
H(\X^*) \leq \polybeta  (\diam (\X^*) + \card (\X^*) )^\polybeta. 
\lbl{polybd}
\eea
Then there exists $\sigma \in (0,\infty)$ such that
 \eq{Hclteq} and \eq{varconv} hold, and $\spann(H) < \infty$ and if 
$b \in (0,\infty)$ with $\spann(H) |b$, 
 then \eq{1205a} holds.
\end{theo}

Now we turn to condition (ii) in Theorem \ref{finranthm}.
Following 
 \cite{PYmfld},
we say that a closed region $A \subset \R^d$ is
a {\em  $d$-dimensional $C^1$ submanifold-with-boundary of
$\R^d$} if it has a differentiable boundary in the following
sense: for every $x$ in the boundary $\partial A$ of $A$,
there is an open  $U \subset \R^d$, and
a continuously differentiable injection $g$ from $U$ to $\R^d$,
such that $0 \in U$ and $g(0)= x$ and
$g(U \cap ([0,\infty) \times \R^{d-1})) = g(U) \cap A$.

\begin{theo}
\lbl{nnlclt}
Let $\kappa \in \N$.  Suppose
$H = H^{(\xi)}$ is induced by a $\xi$ which depends only on the
$\kappa$ nearest neighbours, and
for some $\polybeta  \in (0,\infty)$ suppose we have for all
$(\bx,\X^*)$ that 
\bea
|\xi (\bx;\X^*) | \leq \polybeta (1+R_\kappa(\bx,\X^*))^\polybeta.
\lbl{polynbr}
\eea
 Suppose also that $\supp(f)$ is either a compact convex region 
in $\R^d$ or a compact $d$-dimensional submanifold-with-boundary of
$\R^d$, and suppose $f$ is bounded away from zero on $\supp(f)$.
Finally suppose that the sequence $(r_n)_{n \geq 1}$
satisfies \eq{strongrn}, and that
 $h_{H(\ZZ^*_n)} < \infty$ for some $n \in \N$ with $n > \kappa$. 
Then there exists $\sigma \in (0,\infty)$ such that
\eq{Hclteq} and \eq{varconv} hold,  and
$\spann(H) < \infty$ and if $b \in (0,\infty)$ with
$\spann(H) |b$ then \eq{1205a} also holds.
\end{theo}

We prove Theorems \ref{fin0theo}, \ref{fintheo}
and \ref{nnlclt} in Section \ref{secusePEJP}.
In proving  each of these results,  
we apply  Theorem \ref{finranthm},
and check the CLT condition \eq{Hclteq}
using a general CLT from \cite{PenEJP}, stated below as
 Theorem \ref{lemPEJP}.

  The conclusion that $\sigma >0$ in Theorems \ref{finranthm}--\ref{nnlclt}
and Corollary \ref{finrancoro}
is noteworthy because the result from \cite{PenEJP}
on its own does not guarantee this.
Our approach to showing $\sigma >0$ here is
related to that given in \cite{AB} (and elsewhere)
but is more generic. 
A different approach to providing generic variance lower
bounds was used in \cite{PY1} and \cite{BY} but 
is less well suited to the present setting.

\section{Applications}
\allco

This section contains discussion of some examples of concrete 
models in stochastic geometry, to which the
general local central limit theorems presented
in Section \ref{secstogeo}
 are applicable.
Further
examples 
where the conditions for these
general theorems 
can be verified
are discussed in  \cite{PenEJP,PY1,PY2,PYLLN}. 

\subsection{ Further quantities associated with random geometric graphs}
Suppose the graph $\Gn$ is as in Section \ref{secRGG}.
We assume here that \eq{rhofin} holds with
 $\rho >0$. Theorem \ref{finranthm}  
enables us to  extend the case $\rho >0$ of
 Theorem \ref{Gthm} to non-uniform $f$.
It also yields local CLTs
 for some graph quantities not covered by
 Theorem \ref{Gthm}; we now give some examples.  \\

{\em Number of components for $\Gn$}. This quantity can be written in
the form $H_n(\X_n)$, where $H(\X)$ is the number of components of
the geometric graph $\G(\X,1)$ (which clearly has finite
range interactions).
 In the the thermodynamic limit,
this  quantity satisfies the CLT \eq{0623a} (see Theorem 13.26 of
 \cite{Penbk}). Therefore, 
 Corollary  \ref{finrancoro}
 is applicable here
and shows that it satisfies the local CLT \eq{0623b}. \\

{\em Number of components for $\Gn$ isomorphic to
a given feasible graph $\Gamma$}.
This quantity,
denoted $G^*_n$ in Section \ref{secRGG},
 can be written in the form
$H_n(\X_n)$, with $H(\X)$ the number of
components of $G(\X,1)$ isomorphic to $\Gamma$.
Clearly, this $H$ has finite range interactions since
\eq{finraneq} holds for $\tau =2$.
Also, it satisfies \eq{0623a} by Theorem 3.14 of \cite{Penbk}.
Therefore we can apply Corollary \ref{finrancoro} to deduce \eq{0623b}
 in this case. \\

{\em  Independence number.}
The independence number of a finite graph is the maximal number
$k$ such that there exists a set of $k$ vertices
in the graph   such that none of them are adjacent
Clearly this quantity is 
the sum of the independence numbers of the graph's components,
and therefore if for $\X \subset \R^d$ we set $H(\X)$ to be the independence
 number of ${\cal G}(\X,\tau)$ (also known as the {\em off-line packing
number} since it is the maximum number of balls of radius $\tau/2$
 that can be packed
centred at points of $\X$)  then 
$H$ satisfies  the finite range interactions condition 
\eq{finraneq} with $r=2$.  
Therefore 
we can apply Theorem \ref{fintheo}
to derive a local CLT for the independence number of $\Gn$, as follows.
\begin{theo}.
Let $\tau >0$ and suppose
\eq{subcrit} holds. Suppose $r_n$ is satisfies
\eq{strongrn}. Then 
if
 for $\X \subset \R^d$ we set $H(\X)$ to be the independence
 number of ${\cal G}(\X,\tau)$,
then there exists $\sigma \in (0,\infty)$ such that 
\eq{0623a} holds, and if $b \in \N$ then \eq{0623b} holds.
\end{theo}

\subsection{Germ-grain models}
 Consider a   coverage process in which each point $X_i$ has an associated
mark $T_i$, the  $T_i$ (defined for $i \geq 1$) being i.i.d. nonnegative
random variables with a distribution having
bounded support (i.e., with
$P[T_i \leq K] =1$ for some finite $K$).
Define the random coverage process
\bea
\Xi_n := \cup_{i=1}^n B(r_n^{-1}X_i; T_i).
\lbl{0130a}
\eea
For $U$ a finite union of
convex sets in $\R^d$,
let $|U|$ 
denote the volume
of $U$ (i.e. its Lebesgue measure) and let
 $|\partial U|$
denote the surface area 
of $U$ (i.e. the $(d-1)$-dimensional Hausdorff  
measure of its boundary).

\begin{theo}
Under the above assumptions, if \eq{strongrn} 
holds then
there exists $\sigma >0$ and $\tilde{\sigma} >0$ such that
$n^{-1/2} (|\Xi_n| - \E |\Xi_n|) \tod \NN (0,\sigma^2)$
and 
$n^{-1/2} (|\partial \Xi_n| - \E |\partial \Xi_n|) \tod 
\NN (0,\tilde{\sigma}^2)$, and moreover for any $b \in (0,\infty) $,
\bea
\sup_{u \in \R }
\left\{ \left|
 n^{1/2} P[|\Xi_n|  \in [u,u +b) ] -
\sigma^{-1} b \phi \left(\frac{u - \E |\Xi_n|}{ n^{1/2} \sigma}
\right)
\right|
\right\}
\to 0
~~~~
{\rm as} ~ n \to \infty.
\nonumber \\
\lbl{0623c}
\eea
and
\bea
\sup_{u \in \R }
\left\{ \left|
 n^{1/2} P[|\partial \Xi_n|  \in [u,u +b) ] -
\tilde{\sigma}^{-1} b \phi \left(\frac{u - \E |\partial \Xi_n|}{ n^{1/2} 
\tilde{\sigma}}
\right)
\right|
\right\}
\to 0
~~~~
{\rm as} ~ n \to \infty.
\nonumber \\
\lbl{0623d}
\eea
\end{theo}

{\em Proof.}
The  volume 
 $|\Xi_n|$
can be viewed as a functional $H_n(\X_n^*)$, where
$H(\X) = H^{(\xi)}(\X^*)$ with $\xi((x,t);\X^*)$
given by the volume of that part of the ball
centred at $x$ with radius given by the associated mark $t$,
which is not covered by any corresponding 
ball for some other point $x' \in \X$
with $x'$ preceding $x$ in the lexicographic ordering.
Since we assume
the support of the distribution of the $T_i$ is bounded,
this $\xi$ has finite range $r=2K$.
 Moreover, it satisfies
the polynomial growth bound \eq{polyxi}
so
 by Theorem \ref{fin0theo}
 we get the CLT \eq{Hclteq} and local CLT
\eq{1205a} for any $b >0$ 
(in this example $\spann(H) =0$). 
Thus we have \eq{0623c}.

Turning to the surface area $|\partial \Xi_n|$,
this can also be
viewed as a functional $H_n(\X_n)$ for a different $H = H^{(\xi)}$,
this time taking $\xi(\bx;\X)$ to be the
uncovered surface area of the ball at $x$,
which again has range $r=2K$ and satisfies
\eq{polyxi}.
Hence by Theorem \ref{fin0theo}.
 we get the CLT \eq{Hclteq} and local CLT
\eq{1205a} for any $b >0$ for this choice of $H$ 
(in this example, again $\spann(H) =0$). 
Thus we have \eq{0623d}.
$\qed$ \\

{\em Remark.}
The preceding argument still works if
the independent balls of random radius  in the preceding discussion are
replaced by independent copies of a random compact shape that
is   almost surely contained in the ball $B(K)$ for some $K$
(cf. Section 6.1 of \cite{PenEJP}). \\

{\em Other functionals for the germ-grain model.}
When $f \equiv f_U$, the scaled point process $r_n^{-1/d} \X_n$
can be viewed as a uniform point process in a window of side $r_n^{-1/d}$.
 CLTs for a large class of  other
 functionals on germ-grain models in such a window are considered in
 \cite{HM}, for the Poissonised point process with a Poisson distributed
number of points. Since the Poissonised version of Theorems
\ref{finranthm} and \ref{fin0theo}
 should also hold, it should be possible to
derive local CLTs for many of the quantities considered in \cite{HM},
at least in the case where the  grains (i.e., the balls or other shapes
attached to the random points) are of
uniformly bounded diameter.

\subsection{Random sequential adsorption (RSA).}
RSA (on-line packing) is a model of irreversible deposition of particles
onto an initially empty  $d$-dimensional surface where particles of fixed
finite size arrive sequentially at random locations in
an initially empty region $A$ of a $d$-dimensional
space (typically $d=1$ or $d=2$),
 and each successive particle is  accepted if it does
not overlap any previously  accepted particle.
The region $A$ is taken to be compact and convex.
The locations  of successive particles are independent
and governed by some  density $f$ on $A$. 
In the present setting, we take the mark space $\MM$
to be $[0,1]$ with $\PP_\MM$ 
the uniform distribution. Each
point $\bx = (x,t)$ of $\X^*$
represents an incoming particle with arrival time $t$.
The marks determine the order in which
particles arrive, and two particles at $\bx =(x,t)$ and $\by =(y,u)$
are said to
overlap if $|x-y| \leq 1$.
Let $H(\X^*)$ denote the number of accepted particles.
This choice of
$H$ clearly has finite range interactions (\eq{finraneq} holds for $\tau =2$).

Then $H_n(\X^*_n)$ represents the number of
accepted particles for the re-scaled marked point process $r^{-1}_n \X^*_n$;
note that the density $f$ and hence the region $A$ on which
the particles are deposited, does not vary with $n$.
At least for $r_n = n^{-1/d}$,
 the central limit theorem for $H_n(\X_n)$ 
is known to hold; see \cite{PY2} for the case when 
$A= [0,1]^d$ and $f \equiv f_U$
and \cite{BY} for the extension to the non-uniform case on arbitrary
compact convex $A$ (note that these results do not require the
sub-criticality condition  \eq{subcrit} to be satisfied).
Thus, the $H$ under consideration here satisfies the
condition \eq{Hclteq}.
Therefore we can apply Theorem \ref{finranthm} 
 to obtain a local CLT
for the number of accepted particles in this model.

\begin{theo}
Suppose $f$ has compact convex support and is bounded away
from zero and infinity on its support. 
Suppose $r_n = n^{-1/d}$, and suppose $Z_n = H_n(\X_n^*)$ is the
number of accepted particles in the rescaled RSA model described above.
In other words, suppose
$Z_n$ be the number of accepted particles when 
 RSA is performed on $\X_n$ with distance parameter
$r_n = n^{-1/d}$.
 Then there is a constant $\sigma \in (0,\infty)$   
such that  \eq{normlim1a} holds
and for $b=1$ and $c = n^{1/2}$, \eq{1102c2} holds.
\end{theo}

It is likely that in the preceding result the
condition 
$r_n = n^{-1/d}$ can be relaxed to
\eq{rhofin} holding with $\rho > 0$.
We have not checked the details.

In the {\em infinite input} version of RSA with range of interaction $r$,
 particles continue
to arrive until the region $A$ is saturated, and the total number
of accepted particles is a random variable with its
distribution determined by $r$. 
A central limit theorem  
for the (random) total number of accepted particles
(in the limit  $r \to 0$) is known to hold, 
at least for $f \equiv f_U$; see \cite{SPY}.
It would be interesting to know if a corresponding local central limit
theorem holds here as well.

\subsection{Nearest neighbour functionals}

Many functionals have arisen in the applied literature
which can be expressed as sums of functionals of 
$\kappa$-nearest neighbours,  
for such problems as multidimensional goodness-of-fit tests \cite{BB,BPY},
multidimensional two-sample tests \cite{Henze},
entropy estimation of probability
distributions \cite{LPS}, dimension estimation \cite{LB},
 and nonparametric regression \cite{EJ}.
Functionals considered
 include: sums of power-weighted nearest neighbour distances,
sums of logarithmic functions of the nearest-neighbour distances,
number of nearest-neighbours from the same sample
in a two-sample problem, and others.
Central limit theorems have been obtained explicitly for
some of these examples \cite{BB,Henze,BPY}
and in other cases
 they can often
be derived from more general results
\cite{AB,PenEJP,PY1,Chatt}. Thus, for many of these examples
it should be possible
 to check the conditions of
Theorem \ref{finranthm} (case (ii)).

We consider just one simple example where Theorem \ref{nnlclt}
 is applicable.  Suppose for some fixed $\alpha >0$
that  $H(\X)$ is the
sum of the $\alpha$-power-weighted
nearest neighbour distances in $\X$ (for $\alpha =1$ this
is known as the
total length of the directed nearest neighbour graph on $\X$).
That is, suppose 
$H(\X) = \H^{(\xi)}(\X)$ with $\xi(x;\X)$ given by
 $\min \{|y-x|^\alpha:y \in \X \setminus \{x\}\}$.
 Then $H_n(\X) = r_n^{-\alpha} H(\X)$, and $\xi$ clearly satisfies
\eq{polynbr} for some $\polybeta$, so
 provided $f$ is supported by  a compact convex region
in $\R^d$ or by a compact $d$-dimensional submanifold-with-boundary
of $\R^d$,  and provided $f$ is bounded away from zero on its support,
 Theorem \ref{nnlclt} 
is applicable with $\kappa=1$. Hence in this case
there exists $\sigma \in (0,\infty)$ such that
\eq{Hclteq} and (for any $b \in (0,\infty)$) \eq{1205a} are valid.

\section{Proof of Theorem \ref{genthm2}}
\lbl{secpfgen}
\allco
Let $V, V_1,V_2, V_3,\ldots$ be independent 
identically distributed
 random variables.
Define $\sigma_V := \sqrt{\Var(V)} \in [0,\infty]$.
In the case $\sigma_V=0$, Theorem \ref{genthm2}
is  trivial, so from now on in this section, we assume $\sigma_V > 0$.
Let $b, \mmm_1, \mmm_2, \mmm_3, \ldots$ be positive constants with $\spann_V | b$
and $\mmm_n \sim n^{1/2}$
as $n \to \infty$.

We prove Theorem \ref{genthm2} first
in the special case where $Z_n = S_n$, then
 in the case where $Z_n = Y_n + S_n$,
and then in full generality. Before starting  we recall
a fact about characteristic functions.
\begin{lemm}
\lbl{Vakhlem}
If $\sigma_V = \infty$ then for all $t \in \R$,
as $n \to \infty$
$$
\E \left[ \exp \left(\imag t
n^{-1/2} \sum_{j=1}^n (V_j- \E[V])
\right) \right] \to 0.
$$
\end{lemm}

{\em Proof.} See for example Section 3, and in particular the final
display, of \cite{Vakh}.
 $\qed$
\begin{lemm}
\lbl{classiclem}
Suppose $S_n \eqd \sum_{j=1}^n V_j$ and $\sigma_V < \infty$.
Then as $n \to \infty$,
\bea
\sup_{u \in \R }
\left\{ \left|
\mmm_n P[S_n \in [u,u +b) ] -
\sigma^{-1} b \phi \left(\frac{u - \E S_n}{ \mmm_n
 \sigma_V
}
\right)
\right|
\right\}
\to 0
\lbl{0203a}
\eea
\end{lemm}
{\em Proof.}
First consider
 the special case with $\mmm_n = n^{1/2}$.
In this case, \eq{0203a} holds by
 the classical local central limit theorem
 for sums of i.i.d. non-lattice variables
with finite second moment in the case where $\spann_V=0$
(see
 page 232 of \cite{Breiman}, or Theorem 2.5.4 of
\cite{Durr}),
and by the local central limit theorem for sums of i.i.d. lattice
variables in the case where $\spann_V > 0$ and $b/\spann_V \in \Z$
(see Theorem XV.5.3 of \cite{Feller}, or Theorem 2.5.2 of
\cite{Durr}).

To extend this to the general case with $\mmm_n \sim n^{1/2}$,
observe first that by the special case considered above,
 $n^{1/2} P[S_n \in [u,u+b)]$
remains bounded uniformly in $u$ and $n$, and hence
 \bea
\sup_{u \in \R} \{
| (n^{1/2} - \mmm_n ) P[S_n \in [u,u+b) ] | \}
= \sup_{u \in \R} \left\{
n^{1/2}
\left| 1 - \frac{ \mmm_n }{ n^{1/2} } \right|
P[S_n \in [u,u+b) ] \right\}
\nonumber \\
\to 0. ~~~~~~~~~~
\lbl{1207a}
\eea
Also, for any  $K > 1$,
\bea
 \sup_{|x| \leq K n^{1/2}}
\left\{ \left|
\phi \left(\frac{x}{ n^{1/2} } \right)
- \phi \left(\frac{x}{ \mmm_n } \right)
\right|
\right\}
\leq (2  \pi e)^{-1/2}
 \sup_{|x| \leq K n^{1/2}}
\left\{
\left|
 \left(\frac{x}{ n^{1/2} } \right)
-  \left(\frac{x}{\mmm_n } \right)
\right|
\right\}
\nonumber \\
\leq (2  \pi e)^{-1/2}
\left(
 \frac{K n^{1/2}}{ n^{1/2} }
\right)
\left| 1 -
\frac{ n^{1/2} }{ \mmm_n } \right|
\to 0.
~~~~~~~~~~~
\lbl{1207b}
\eea
Also, for large enough $n$,
\bean
\sup_{|x| \geq K n^{1/2}} \max
\left( \phi \left( \frac{x}{\mmm_n }  \right),
 \phi \left( \frac{x}{n^{1/2} }  \right) \right) \leq \phi(K-1)
\eean
and since $K$ is arbitrarily large,
 combined with \eq{1207b}, this shows that
\bean
 \sup_{x  \in \R }
\left\{ \left| \phi \left(\frac{x}{ n^{1/2} } \right)
- \phi \left(\frac{x}{ \mmm_n } \right) \right| \right\} \to 0.
\eean
Combined with \eq{1207a}, this shows that we can deduce
\eq{0203a}
for general $\mmm_n $ satisfying $\mmm_n \sim n^{1/2}$
from the special case with $\mmm_n = n^{1/2}$ which was
established earlier.  $\qed$ \\

\begin{lemm}
\lbl{lemZXY}
Theorem \ref{genthm2} holds in the special case where $Z_n = Y_n + S_n$.
\end{lemm}
{\em Proof.}
Assume, along with the hypotheses of Theorem \ref{genthm2}, that
$Z_n = Y_n + S_n$.  Considering characteristic functions,
by \eq{normlim1a} we have for $t \in \R$ that
\bea
\E \left[ \exp \left(
\imag t n^{-1/2} (Y_n- \E Y_n)
 \right) \right] \E \left[ \exp \left(
 i t n^{-1/2} (S_n- \E S_n) \right) \right]
\nonumber \\
\to \exp ( - \frac{1}{2} t \sigma^2).
~~~~~~~
\label{cfeq}
\eea
If $\sigma_V  = \infty$ then
 by Lemma \ref{Vakhlem}, the second factor
in the  left hand side of \eq{cfeq}
tends to zero, giving a contradiction. Hence we may assume $\sigma_V < \infty$
from now on.

 By the Central Limit Theorem,
\bea
n^{-1/2} ( S_n - \E S_n )
 \tod N(0, \sigma_V^2).
\lbl{normlim2a}
\eea
By \eq{cfeq} and \eq{normlim2a}, $\sigma_V^2 \leq \sigma^2$ and
 setting $\sigma_Y^2 := \sigma^2 - \sigma_V^2 \geq 0$, we have that
$ n^{-1/2} (Y_n- \E Y_n) $
is asymptotically $ \NN(0,\sigma_Y^2 ).  $ Hence,
\bea
\mmm_n^{-1} (Y_n- \E Y_n) \tod \NN(0,\sigma_Y^2 ).
\lbl{normlimX}
\eea
 That is, \eq{0110a} holds.

Let $u \in \R$ and set
\bea
t := t(u,n) := \mmm_n^{-1}(u - \E Z_n).
\lbl{0204a}
\eea
Assume 
 that $Z_n = Y_n + S_n$.
By independence of $Y_n$ and $S_n$,
\bean
 P[Z_n \in [u,u +b) ]
=
 P [ \mmm_n^{-1} ( Z_n  - \E[Z_n])  \in \mmm_n^{-1}[u-\E Z_n,u+b - \E Z_n) ]
\\
= \int_{-\infty}^{\infty} P \left[ \frac{Y_n - \E Y_n  }{\mmm_n} \in dx
\right]
 P\left[ \frac{S_n - \E S_n }{\mmm_n} \in
 \mmm_n^{-1}[u-\E Z_n,u+b - \E Z_n) -x
\right]
\eean
so that
\bean
 \mmm_n
 P[Z_n \in [u,u +b) ]
= \int_{-\infty}^{\infty} P \left[ \frac{Y_n - \E Y_n  }{\mmm_n} \in dx
\right]
\\
\times
\left(  \mmm_n P\left[ S_n - \E S_n \in
 [u-\E Z_n - x \mmm_n,u - \E Z_n -x  \mmm_n +b)
\right]
\right)
\\
= \int_{-\infty}^{\infty} P \left[ \frac{Y_n - \E Y_n  }{ \mmm_n} \in dx
\right]
\left(   \mmm_n P\left[ S_n - \E S_n \in
[ (t - x )  \mmm_n,(t -x)  \mmm_n +b)
\right]
\right) .
\eean
By Lemma \ref{classiclem},
\bean
 \mmm_n P\left[ S_n - \E S_n \in [ y  \mmm_n, y  \mmm_n +b) \right] =
 \frac{b}{\sigma_V} \phi \left( \frac{ y}{\sigma_V} \right) + g_n(y)
\eean
where
\bea
\sup_{y \in \R} | g_n(y) | \to 0 ~~~~~ {\rm as}~~~ n \to \infty.
 \lbl{1102b}
\eea
Hence,
\bean
  \mmm_n P[Z_n \in [u,u +b) ]  = \E \left[
 \frac{b}{\sigma_V} \phi \left( \frac{ t-  \mmm_n^{-1} (Y_n-\E Y_n) }{
\sigma_V } \right) +
 g_n \left( t-  \mmm_n^{-1}(Y_n - \E Y_n) \right)
\right],
\eean
so by \eq{1102b}, to prove \eq{1102c2},
 it suffices to prove
\bea
\sup_{u \in \R }
\left\{
\left|
 \E \left[ \sigma_V^{-1} \phi \left( \frac{t(u,n)-
 \mmm_n^{-1} (Y_n - \E Y_n )}{ \sigma_V} \right) \right]
-
\sigma^{-1}\phi \left(\frac{u - \E Z_n}{  \mmm_n \sigma}
\right)
\right|
\right\}
\to 0. \nonumber \\
\lbl{1102d}
\eea
Suppose this fails. Then there is a strictly increasing sequence of
natural numbers $(n(m), m \geq 1)$
 and  a
 sequence of real numbers $(u_m,m \geq 1)$
 such that with $t_m := t(u_m,n(m)),$ we have
\bea
\liminf_{m \to \infty}
\left|
 \E \left[ \sigma_V^{-1} \phi \left( \frac{t_m -
 \mmm_{n(m)}^{-1} (Y_{n(m)} - \E Y_{n(m)} )}{ \sigma_V} \right) \right]
-
\sigma^{-1}\phi \left(\frac{u_m - \E Z_{n(m)}}{  \mmm_{n(m)} \sigma}
\right)
\right|
> 0.
\nonumber \\
\lbl{1102e}
\eea

By taking a subsequence if necessary,
we may assume without loss of generality, either that
$t_m \to t $ for some $t \in \R$, or that
or that $|t_m | \to \infty$ as $m \to \infty$.
Consider first the latter case. If
 $|t_m | \to \infty$ as $m \to \infty$, then
by \eq{normlimX},
\bean
 P[
|t_m -   \mmm_{n(m)}^{-1}(Y_{n(m)}- \E Y_{n(m)}) | \leq  |t_m |/2]
\leq P [ |   \mmm_{n(m)}^{-1}(Y_{n(m)}- \E Y_{n(m)}) | \geq |t_m |/2]
\\
\to 0,
\eean
and hence
$$
 \E \left[ \sigma_V^{-1} \phi \left( \frac{t_m -
 \mmm_{n(m)}^{-1} (Y_{n(m)} - \E Y_{n(m)} )}{ \sigma_V} \right) \right]
\to 0.
$$
Since $ \mmm_{n(m)}^{-1}(u_m - \E Z_{n(m) })$ is
equal to $t_{m}$ by \eq{0204a}, we also have under this assumption
that
$\sigma^{-1}\phi \left(\frac{u_m - \E Z_{n(m)}}{  \mmm_{n(m)} \sigma}
\right)$ tends to zero, and thus we obtain
a contradiction of \eq{1102e}.

In the  case where $t_m \to t$ for some finite $t$, we have
by \eq{normlimX} that
$ t_m -    \mmm_{n(m)}^{-1}(Y_{n(m)}- \E Y_{n(m)}) $ converges in distribution
to $t-W_1$, where
$W_1 \sim \NN(0, \sigma_Y^2 )$. Hence
as $m \to \infty$,
\bean
 \E \left[ \sigma_V^{-1} \phi \left( \frac{t_m -  \mmm_{n(m)}^{-1}(Y_{n(m)} - \E
Y_{n(m)} )  }{ \sigma_V } \right)
\right]
\to \sigma_V^{-1} \E \phi((t-W_1)/ \sigma_V)
\\
= \E f_{W_2} (t - W_1),
\eean
where $W_2 \sim N(0,\sigma_V^2)$, with
 probability density function
 $f_{W_2}(x) : = \sigma_V^{-1} \phi(x/\sigma_V)$.
If we
 assume $W_1$, $W_2$ are independent,
then $\E f_{W_2} (t - W_1)$ is the convolution formula for the
probability density function of $W_1 + W_2$,
which is
 $ \NN(0, \sigma^2)$, so that
$$
 \E f_{W_2} (t - W_1)
= f_{W_1 + W_2}(t) = \sigma^{-1}\phi(t/\sigma).
$$
On the other hand, since $ \mmm_{n(m)}^{-1}(u_m - \E Z_{n(m) })$ is
equal  (by \eq{0204a}) to $t_{m}$ which we assume converges to $t$,
we also have that
$$
\sigma^{-1}\phi \left(\frac{u_m - \E Z_{n(m)}}{  \mmm_{n(m)} \sigma}
\right)
\to
\sigma^{-1}\phi \left(\frac{t}{ \sigma}
\right),
$$
and therefore we obtain a contradiction of
\eq{1102e} in this case too.

 Thus
\eq{1102e} fails, and therefore
\eq{1102d} holds. Hence,
 \eq{1102c2} holds in the case with
 $Z_n = Y_n + S_n$.
$\qed$ \\

{\em Proof of Theorem \ref{genthm2}.}
Set $Z'_n:=Y_n + S_n$.  By the integrability assumptions,
$Z'_n$ is integrable. By
\eq{normlim1a} and the assumption that
$n^{-1/2}\E[|Z_n-Z'_n|] \to 0$ as $n \to \infty$,
\bea
n^{-1/2} (Z'_n- \E Z'_n) \tod \NN(0,\sigma^2 )
{\rm ~~~as ~~} n \to \infty.
\eea
Let $b >0$ with $\spann_V | b$.  By Lemma \ref{lemZXY},
$\sigma^2 \geq \Var V$ and \eq{0110a} holds and
\bean
\sup_{u \in \R }
\left\{ \left|
 \mmm_n
 P[Z'_n \in [u,u +b) ] -
\sigma^{-1} b \phi \left(\frac{u - \E Z'_n}{  \mmm_n
 \sigma
}
\right)
\right|
\right\}
\to 0
~~~~
{\rm ~~~as ~~} n \to \infty.
\eean
Hence, by the assumption $n^{1/2} P[Z_n \neq Z'_n] \to 0$,
\bean
\sup_{u \in \R }
\left\{ \left|
 \mmm_n
 P[Z_n \in [u,u +b) ] -
\sigma^{-1} b \phi \left(\frac{u - \E Z'_n}{  \mmm_n
 \sigma
}
\right)
\right|
\right\}
\to 0
~~~~
{\rm ~~~as ~~} n \to \infty,
\eean
and since the assumption
$n^{-1/2}\E[|Z_n-Z'_n|] \to 0$ implies that
$ \mmm_n^{-1}(\E[Z_n] - \E [Z'_n]) \to 0$
as $n \to \infty$, and $\phi$ is uniformly continuous on $\R$,
we can then deduce \eq{1102c2}.
 $\qed$

\section{Proof of theorems for percolation}
\lbl{secpfperc}
\allco
We shall repeatedly use the following Chernoff-type tail bounds
for the binomial and Poisson distributions 
For $a >0$ set $\varphi(a):=1-a + a \log a$.
Then $\varphi(1)=0$ and $\varphi(a) >0$ 
for $a \in (0,\infty)\setminus \{1\}$.
\begin{lemm}
\lbl{PenL1.1}
If  $X$ is a binomial or Poisson distributed random variable
with  $\E[X] = \mu >0$.
Then we have for all $x > 0$ that
\bea
P[X \geq x] \leq \exp( - \mu \varphi(x/\mu)), ~~~~
x \geq \mu; \lbl{Biup}
\\
P[X \leq x] \leq \exp( - \mu \varphi(x/\mu)),  ~~~~
x \leq \mu. \lbl{Bilo}
\eea
\end{lemm}
{\em Proof.} See e.g. Lemmas 1.1 and 1.2  of \cite{Penbk}. $\qed$ \\

{\em Proof of Theorem \ref{thclus}.}
Let $(B_n)_{n \geq 1}$ be a sequence of nonempty finite subsets
in $\Z^d$ with vanishing relative boundary. The first conclusion
\eq{ClusCLT} follows from Theorem 3.1 of \cite{PenCLT}, so it
remains to prove \eq{ClusLLT}.

For $x \in \Z^d$ let $\|x\|_\infty$ denote the $\ell_\infty$-norm of $x$,
i.e., the maximum absolute value of its coordinates.
Let $B_n^{o}$
be the set of points $x$ in $B_n$ such that all $y\in\Z^d$  with
$\|y-x\|_\infty \leq 1$ are also in $B_n$.
Since $|B_n \setminus B_n^o|/|\partial B_n|$ is bounded by a constant
depending only on $d$, the vanishing
relative boundary condition
\eq{vrb} implies $|B^o_n|/|B_n| \to 1$ as $n \to \infty$.

 Hence, by the pigeonhole principle,
 for all large enough $n$ we can choose a set of
points $x_{n,1}, x_{n,2}, \ldots, x_{n,\lfloor 5^{-d} |B_n|/2\rfloor}$
in $B^o_n$ such that  $\|x_{n,j} - x_{n,k}\|_\infty \geq 3$
for each distinct $j,k$ in $\{1,2,\ldots,\lfloor  5^{-d} |B_n|/2\rfloor\}$
(let these points  be chosen by some arbitrary deterministic rule).

For $1 \leq j \leq \lfloor 5^{-d} |B_n|/2\rfloor$, let $I_{n,j}$
be the indicator of the event that each vertex $y \in \Z^d$ with
$\|y-x_{n,j}\|_\infty =1$ is closed,
 and list the $j$ for which $I_{n,j} =1$, in increasing order,
as $J(n,1)  \ldots, J(n,N_n)$, where
 $N_{n} := \sum_{j=1}^{ \lfloor 5^{-d} |B_n|/2\rfloor} I_{n,j}$.
Let $I'_{n,j}$ be the indicator of the event that the vertex
 $x_{n,j}$ is itself open.
 Then $N_n$ is binomially distributed with
parameter
$(1-p)^{3^d -1}$,
 so by Lemma \ref{PenL1.1},
\bea
\lbl{0129a}
\limsup_{n \to \infty} |B_n|^{-1} \log  P[
N_n <
5^{-d}  (1-p)^{3^d -1} |B_n|/4]  <0.
\eea
Set $b_n:=
\lfloor 5^{-d}  (1-p)^{3^d -1} |B_n|/4 \rfloor.
$
Let $V_1,V_2,\ldots$ be a sequence of independent Bernoulli
variables with parameter $p$, independent of
everything else.
 Recalling that $\Cl(B)$ denotes the  number of open clusters in $B$,
set
\bean
S'_n :=
\sum_{j=1}^{\min(b_n ,N_n)
}
I'_{n,J(n,j)};  ~~~~~
Y_n:= \Cl(B_n)-S'_n,
\eean
and
\bean
S_n := S'_n + \sum_{j=1}^{(b_n - N_n)^+} V_j,
\eean
where $x^+ : = \max(x,0)$ as usual, and the sum $\sum_{i=1}^0$ is
taken to be zero.

In this case, the `good boxes' discussed in Section \ref{secintro}
are the unit $\ell_\infty$-neighbourhoods of the sites
$x_{n,J(n,1)},
x_{n,J(n,2)}, \ldots
x_{n,J(n,\min(b_n,N_n))}
$. If $x_{n,j}$ is at the centre of a good box, it is (if open)
isolated from other open sites, so that $Y_n$ is simply the number
of open clusters in $B_n$  if one ignores
 all sites $x_{n,J(n,j)}$ $(1 \leq j \leq
\min(b_n,N_n))$.
 Hence $Y_n$ does not affect
the open/closed status of these sites.

Thus $S_n$ has the $\Bin(b_n,p)$ distribution and its distribution,
given $Y_n$, is unaffected by the value of $Y_n$ so $S_n$
is independent of $Y_n$.
Also,
 \bean
\Cl(B_n) - ( Y_n + S_n) = S'_n -S_n = - \sum_{j=1}^{(b_n-N_n)^+} V_j
\eean
so that by \eq{0129a}, both $|B_n|^{1/2} P[\Cl(B_n) \neq Y_n+S_n]$
and $|B_n|^{-1/2}\E[|\Cl(B_n) -(Y_n+S_n)|]$ tend to
zero as $n \to \infty$. Combined with \eq{ClusCLT} this shows
that Theorem \ref{genthm2} is applicable, with $\spann_V =1$,
and that result shows that
\eq{ClusLLT} holds. $\qed$ \\

In the proof of Theorem \ref{thlargest}, and again later  on,
we shall use the following.
\begin{lemm}
\lbl{Azulem}
Suppose $\xi_1,\ldots,\xi_m$ are independent identically
distributed random elements
of some measurable space $(E,{\cal E})$.
Suppose $m \in \N$ and $\psi:E^m \to \R$ is measurable and
suppose for some finite $K$ that for $j =1,\ldots,m$,
$$
K \geq \sup_{(x_1,\ldots,x_m,x'_j) \in E^{m+1}}
|\psi(x_1,\ldots,x_j,\ldots,x_m)- \psi(x_1,\ldots,x'_j,\ldots,x_m)|.
$$
Set $Y = \psi (\xi_1,\ldots,\xi_m)$.
Then
for any $t >0$,
$$
P[ |Y- \E Y| \geq t ] \leq 2 \exp ( - t^2/(2 m K^2)).
$$
\end{lemm}
{\em Proof.} 
The argument is similar to e.g. the proof of Theorem 3.15 of \cite{Penbk};
 we include it for completeness.  For $1 \leq i \leq m$
let $\F_i$ be the $\sigma$-algebra generated
by $\xi_1,\ldots, \xi_i$, and let $\F_0$ be the
trivial $\sigma$-algebra. Then $Y- \E[Y]= \sum_{i=1}^m D_i$
with $D_i := \E[Y|\F_i] - \E[Y|\F_{i-1}]$, the $i$th martingale
difference. Then with $\xi'_i$ independent of $\xi_1,\ldots,\xi_m$
with the same distribution  as them, we have
$$
D_i = \E[\psi(\xi_1,\ldots, \xi_i, \ldots \xi_m) -
 \xi(\xi_1,\ldots,\xi'_i, \ldots, \xi_m)|\F_i] 
$$  
so that  $|D_i| \leq K$ almost surely and hence by 
Azuma's inequality (see e.g. \cite{Penbk}) we
have the result.  \\

{\em Proof of Theorem \ref{thlargest}.}
Assume $d \geq 2 $ and $p > p_c(d)$.
Let $(B_n)_{n \geq 1}$ be a cube-like sequence of lattice boxes in $\Z^d$.
For finite nonempty $A \subset \Z^d$ we
define the {\em diameter} of $A$, written $\diam(A)$, to be
  $\max \{\|x-y\|_\infty: x \in A, y \in A\}$.

Set $\bxdi_n:= \lceil \diam(B_n)^{1/(4d)}\rceil$.
Let $B_n^{{\rm in}}$ be the set of points
 $x$ in $B_n$ such that all $y\in\Z^d$  with
$\|y-x\|_\infty \leq \bxdi_n $ are also in $B_n$.
Then we claim that
 $|B^{{\rm in}}_n|/|B_n| \to 1$ as $n \to \infty$.
Indeed, writing  $B_n =  \prod_{j=1}^d([-a_{j,n},b_{j,n} ] \cap \Z)$, 
 from the cube-like condition \eq{cubelike} we have for
$1 \leq j \leq d$ that
$\bxdi_n = o(a_{j,n} + b_{j,n})$ as $n \to \infty$, and therefore
\bean
 |B^{{\rm in}}_n| =
 \prod_{j=1}^d (b_{j,n} +a_{j,n} - 2 \bxdi_n)   =
 (1 +o(1)) 
 \prod_{j=1}^d (a_{j,n} +b_{j,n} ), 
\eean
justifying the claim.

By the preceding claim, and the pigeonhole principle,
 for all large enough $n$ there is a deterministic set of points
$x_{n,1}, x_{n,2}, \ldots, x_{n,\lfloor 5^{-d} |B_n|/2\rfloor}$
in $B^{{\rm in}}_n$ such that  $\|x_{n,j} - x_{n,k}\|_\infty \geq 3$
for each distinct $j,k$ in $\{1,2,\ldots,\lfloor  5^{-d} |B_n|/2\rfloor\}$.

For $1 \leq j \leq \lfloor 5^{-d} |B_n|/2 \rfloor$, let $I_{n,j}$
be the indicator of the event that (i) each vertex
$y \in \Z^d$ with $\|y-x_{n,j}\|_\infty =1$ is open,    and (ii)
the open cluster in $B_n$ containing all
$y \in \Z^d$ with $\|y-x_{n,j}\|_\infty =1$
has diameter at least $\bxdi_n$.

Set $m(n) : = \lfloor 5^{-d} p^{3^d-1} \theta_d(p) |B_n|/8 \rfloor$,
 with $\theta_d(p)$ denoting the percolation probability.
List the $j$ for which $I_{n,j} =1$
as $J(n,1),$ $ \ldots,$ $J(n,N_n)$, with
 $N_{n} := \sum_{j=1}^{ \lfloor 5^{-d} |B_n|/2\rfloor} I_{n,j}.$
  Then
we have for $n$ large that
$$
\E [ N_n ] \geq
\lfloor 5^{-d} |B_n|/2\rfloor
p^{3^d-1} \theta_d(p) \geq 2 m(n).
$$
Changing the open/closed status of a single site
$z$ in $B_n$ can change the value of $I_{n,j}$  only for
 those $j$ for which
$\|x_{n,j} -z\|_\infty \leq \bxdi_n$, and the number of such $j$
is at most $(2 \bxdi_n +1)^d$. Moreover,
for $n$ large
$$
(2 \bxdi_n +1)^d \leq ( 2 (\diam B_n)^{1/(4d)} +3 )^d \leq
3^d (\diam B_n)^{1/4}
\leq 3^d |B_n|^{1/4}
$$
 so that the total change in $N_n$ due to
changing the status of a single site $z$ is at most
 $3^d|B_n|^{1/4}$.
So by Lemma \ref{Azulem},
\bean
P[ N_n \leq
m(n)
 ]  
\leq
P [ |N_n- \E N_n| \geq m(n)
]
\leq 2 \exp  \left( - \frac{ m(n)^2
 }{  2 |B_n| (3^d|B_n|^{1/4})^2 }
\right)
\eean
and hence
\bea
\lbl{0129b}
\limsup_{n \to \infty} |B_n|^{-1/2} \log  P[
N_n \leq
m(n)
]  <0.
\eea
Let $V_1,V_2,\ldots$ be a sequence of independent Bernoulli
variables with parameter $p$, independent of everything else.
For $1 \leq j \leq \lfloor 5^{-d} |B_n|/2\rfloor$, let 
 $I'_{n,j}$ be the indicator of the event that the vertex $x_{n,j}$
 is open.
  Set
\bean
S'_n :=
\sum_{j=1}^{\min(m(n) ,N_n)
}
I'_{n,J(n,j)}; ~~~~~
S_n := S'_n + \sum_{j=1}^{(m(n) - N_n)^+} V_j.
\eean

Let $Y_n$ be the size of the largest open cluster in $B_n$ if
the status of  $x_{i,n}$ is set to `closed' for the first
$\min(m(n),N_n)$ values of $j$ for which $I_{n,j} =1$.

Then $S_n$ has the $\Bin(m(n),p)$ distribution and we assert that
its distribution,
given $Y_n$, is unaffected by the value of $Y_n$ so $S_n$
is independent of $Y_n$. Indeed, $Y_n$ is obtained without sampling
the status of the  sites $x_{n,j}$ for
the first
$\min(m(n),N_n)$ values of $j$ for which $I_{n,j} =1$.

To explain this further, consider algorithmically sampling the
open/closed status
of sites in $B_n$ as follows.
 First sample the status of
sites outside $\cup_{j}\{x_{n,j}\}$. Then sample the status
of those $x_{n,j}$ for which the $\ell_\infty$-neighbouring sites are not
all open (for these sites, $I_{n,j}$ must be zero).
At this stage,
it remains to sample the status of sites
$x_{n,j}$ for which the $\ell_\infty$-neighbouring sites are all open,
and for these sites one can tell, without revealing the
value of $x_{n,j}$, whether or not $I_{n,j} =1$ (and in particular
one can determine the value of $N_{n}$).
At the next step sample the status of all $x_{n,i}$
except for the first $\min(N_{n},m(n))$
values of $i$ which have $I_{n,j} =1$. At this point, the value
of $Y_n$ is determined. However, the value of $S_n$ is
determined by the status of the remaining unsampled sites
together with some extra Bernoulli variables in the
case where $N_{n} < m(n)$, so its distribution is independent
of the value of $Y_n$ as asserted.

Next, we establish that $\LL(B_n) = Y_n +S_n$
with high probability.
One way in which this could fail would be if $N_n < m(n)$,
but we know from \eq{0129b} that this has small probability.
Also, we claim that with high probability, all sites  $x_{n,j}$
for which $I_{n,j} =1$ have all their neighbouring sites
as part of the largest open cluster, regardless
of the status of $x_{n,i}$. To see this, let $A_n$ be the event
that (i) there is a unique open cluster for $B_n$
that crosses $B_n$ in all directions (in the
sense of \cite{PenCLT}) and
(ii) all other clusters in $B_n$ have diameter less than
$\bxdi_n$. Then we claim that
 $P[A^c_n]$ decays exponentially in $\gamma_n$ in the sense that
\bea
\limsup_{n \to \infty} (\diam B_n)^{1/(4d)} \log P[A_n^c] < 0.
\lbl{0129c}
\eea
 The proof of \eq{0129c} proceeds as in  
proof of Lemma 3.4 of \cite{PenCLT}; we include a sketch
of this argument here for completeness.

First suppose $d=2$.  For a given rectangle of
 dimensions $(\gamma_n/3) \times \gamma_n$,
the probability that it fails to have an open  crossing the long way 
decays exponentially in $\gamma_n$ (see Lemma 3.1 of \cite{PenCLT}).
Consider the family of all rectangles 
of dimensions $(\gamma_n/3) \times \gamma_n$ or of dimensions
 $\gamma_n \times (\gamma_n/3)$,
with all corners in $(\gamma_n/3) \Z^2$, having non-empty intersection
with $B_n$. The number of such rectangles is $O(\diam(B_n)^{d-1/2})$.   
By the preceding probability estimate,
all  rectangles in this family
 have an open crossing
the  long way, except on an event 
of probablity decaying exponentially in $\gamma_n$. However,
if all these rectangles have an open crossing the long way, then
event $A_n$ occurs and we have justified \eq{0129c} 
for $d=2$.

For $d\geq3$, by the well known result of Grimmett and Marstrand
\cite{GM},
there exists a finite $K$
such that there is an infinite
open cluster in the slab $[0,K] \times \R^{d-1}$ with
 strictly positive probability.
By dividing $B_n$ into slabs of thickness
$K$ we see for $1 \leq i \leq d$
 that the probabilty that there is no open  crossing of $B_n$
 in the $i$-direction decays exponentially
in $\diam(B_n)$.  Moreover, for $i \neq j$,
by a similar slab argument
(consider successive slabs of thickness $K$ in the $i$ direction),
the probability that there is an open cluster
in $B_n$ that  crosses $B_n$ in the $i$ direction
but not the $j$ direction decays exponentially in $\diam(B_n)$.
Similarly the probability that there are two or more disjoint
open clusters in $B_n$ which  cross in the $i$ direction
 decays exponentially in $n$.
Finally by a further slab argument,
the probability that there is an open cluster 
which has diameter at least $\gamma_n/d$ in the $i$ direction
  but fails to cross the whole of $B_n$ in the $j$ direction,
 decreases exponentially in $\gamma_n$.  This 
 justifies \eq{0129c} for $d\geq 3$.

Note that the occurrence or otherwise of $A_n$
is unaffected by the open/closed status of those $x_{n,i}$ for which
$I_{n,j}=1$.
Also, for large enough $n$,  on
 event $A_n$, whatever status we give to
 these $x_{n,j}$, the unique crossing cluster is the largest one
because  it has at least $\diam(B_n)$ elements while all other
clusters have at most $O(\diam(B_n)^{1/4})$ elements.

If $N_n \geq m(n)$ and event $A_n$ occurs, then for each $j\leq m(n)$,
the site $x_{n,J(n,j)}$ is in the largest open cluster
if and only if it is open, since if it is open then it is in
an open cluster of diameter at least $\bxdi_n$. This shows that
 if $N_n \geq m(n)$ and event $A_n$ occurs,
we do indeed have
 $\LL(B_n) = Y_n + S_n$.
 Together with
 the previous probability
estimates \eq{0129b} and \eq{0129c}, this shows that
$|B_n|^{1/2} P [\LL(B_n) \neq Y_n + S_n] \to 0$ as $n\to \infty$.
Moreover, by the Cauchy-Schwarz inequality,
\bean
 \E [|\LL(B_n) - ( Y_n + S_n)|]
=
 \E [|\LL(B_n) - ( Y_n + S_n)|{\bf 1}_{\{N_n < m(n) \} \cup A_n^c}]
\\
\leq ( P[N_n < m(n)] + P[A_n^c] )^{1/2} (\E[ (\LL(B_n) -( Y_n + S_n))^2])^{1/2}
\\
\leq ( P[N_n < m(n)] + P[A_n^c] )^{1/2} (|B_n| + m(n)) \to 0.
\eean
By Theorem 3.2 of \cite{PenCLT}, the first conclusion \eq{LargCLT}
holds,  and by the preceding discussion, we can then apply
Theorem \ref{genthm2} with $\spann_V=1$,
 to derive the second conclusion \eq{LargLLT}.
 $\qed$

\section{Proof of Theorem \ref{Gthm}}
\lbl{secpfedges}
\allco
We are now in the setting of Section \ref{secRGG}.
Assume $f \equiv f_U$, and fix a feasible connected graph $\Gamma$ with
$\kappa$ vertices $(2 \leq \kappa < \infty)$.
 Assume also that the sequence $(r_n)_{n \geq 1}$ is given
and satisifies \eq{rhofin} and \eq{taubig}.
Then
 $P[\G(\X_\kappa,1/(\kappa+3)) \sim \Gamma] \in ( 0,1)$.
Let $Q_{n,1},Q_{n,2},\ldots,Q_{n,m(n)}$ be disjoint cubes
 of side $(\kappa+5)r_n$, contained in the unit cube,
with $m(n) \sim ((\kappa+5) r_n)^{-d}$ as $n \to  \infty$.
For $1 \leq j \leq m(n)$, let
$I_{n,j}$ be the indicator of the event
that
$\X_n  \cap Q_{n,j}$
consists of exactly $\kappa$ points, all of them at a
 Euclidean distance greater than $r_n$ from the boundary of $Q_{n,j}$.
List the indices $j \leq m(n)$ such that
$I_{n,j}=1$,
in increasing order, as
$J_{n,1},\ldots,J_{n,N_n}$, with $N_n: = \sum_{j=1}^{m(n)} I_{n,j}$.
Then
\bea
\E [N_n] = m(n)  ((\kappa+3)/(\kappa+5))^{d\kappa}
P[\Bin(n, ((\kappa+5)r_n)^d) =\kappa] ,
\lbl{1027a}
\eea
and hence as $n \to \infty$,  since
$n r_n^d$ is bounded by our assumption \eq{rhofin},
\bea
\E [N_n]
 \sim
\kappa!^{-1}(\kappa+3)^{d\kappa}
(\kappa+5)^{-d}
n^\kappa r_n^{d(\kappa-1)}
 \exp(-n (\kappa+5)^d r_n^d).
\lbl{1113a}
\eea
 Recalling from \eq{taubig} that $\tau_n := \sqrt{ n (nr_n^d)^{\kappa-1}}$,
we can rewrite \eq{1113a} as
\bea
\E [N_n]
\sim  
\kappa!^{-1}(\kappa+3)^{d\kappa}
(\kappa+5)^{-d}
\tau_n^2
 \exp(-n (\kappa+5)^d r_n^d)
\lbl{1031a}
\eea
as $n \to \infty$.
Moreover, for the Poissonised version of this model
where the number of points is Poisson distributed
with mean $n$, we have  the same asymptotics
for the quantity corresponding to $N_n$ (the binomial
probability in \eq{1027a}  is asymptotic to the corresponding
Poisson probability).
Set  $\alpha$ to be one-quarter of the coefficient of
$\tau_n^2$ in \eq{1031a}, if the exponential factor is
replaced by its smallest value in the sequence, i.e. set
\bea
\alpha :=
(4 \kappa!)^{-1}(\kappa+3)^{d\kappa}
(\kappa+5)^{-d}
\inf_n \exp( -n (\kappa+5)^d r_n^d).
\lbl{0113c}
\eea
Then
 $\alpha > 0$ by our assumption \eq{rhofin} on $r_n$.
\begin{lemm}
\lbl{lem1}
It is the case that
$$
\limsup_{n \to \infty} \tau_n^{-2} \log
P \left[ N_n < \alpha \tau_n^2 \right]   < 0.
$$
\end{lemm}
{\em Proof.}
Let $\delta >0$ (to be chosen later).
Let $M_n$ be Poisson distributed with parameter
$(1-\delta)n$, independent of the sequence
of random $d$-vectors
$X_1,X_2,\ldots$. Define the 
Poisson point process 
$$
 \Po_{n(1-\delta)} := \{X_1,\ldots,X_{M_n}\} .
$$
 Let $N'_n $ be defined in the same manner as $N_n$
but in terms of $\Po_{n(1-\delta)}$ rather than $\X_n$.
That is, set
$$
N'_n := \sum_{j=1}^{m(n)} I'_{n,j}
$$
with $I'_{n,j}$ denoting the indicator
of the event that $\Po_{n(1-\delta)} \cap Q_{n,j}$
consists of exactly $\kappa$ points, all
at distance greater than $r_n$ from the boundary
of $Q_{n,j}$.
List the indices $j \leq M_n$ such that
 $I'_{n,j}=1$ as $J'_{n,1}, \ldots, J'_{n,N'_n}$.

Since \eq{1031a} holds in the Poisson setting too,
using the definition of $\tau_n$ we have as $n \to \infty $ that
\bea
\E [N'_n]
\sim  
\kappa!^{-1}(\kappa+3)^{d\kappa}
(\kappa+5)^{-d}
(1-\delta)^\kappa
\tau_n^2
 \exp(-n (1-\delta) (\kappa+5)^d r_n^d).
\lbl{0113b}
\eea
 By \eq{1031a} and \eq{0113b},
  we can and do choose $\delta >0$ to be small enough
so that
 $\E N'_n > (3/4) \E N_n$
 for large $n$. 

 By \eq{1031a} and \eq{0113c} we have for large $n$
that $2 \alpha \tau_n^2 \leq (5/8) \E N_n$.
Also, $N'_n$ is binomially distributed,
and hence by Lemma \ref{PenL1.1},
$P[N'_n < 2 \alpha \tau_n^2 ]$
decays exponentially in  $\tau_n^2$.

By Lemma \ref{PenL1.1},
except on an event of probability
decaying exponentially in  $n$, the value of $M_n$ lies between
$n(1 - 2 \delta)$
and
$n$.
If this happens, the discrepancy between $N_n$ and $N'_n$
is due to the addition of at most an extra $2 \delta n$ points
to $\Po_{n(1-\delta)}$.  
If also
$N'_n \geq 2 \alpha \tau_n^2$ then to have
$N_n < \alpha \tau_n^2$, at least $\alpha \tau_n^2$
of the added points must land in the union of the first
$ \lceil 2 \alpha \tau_n^2 \rceil$ cubes
 contributing to $N'_n$. 

To spell out the preceding argument in more detail,
 let $1 \leq j \leq m(n)$.
If $M_n < n$ and $I'_{n,j}=1$ and $X_k \notin Q_{n,j}$
for $M_n < k \leq n$, then $I_{n,j}=1$, since in
this case $
\X_{n} \cap Q_{n,j}
 = 
\Po_{n(1-\delta)} \cap Q_{n,j}
 $. 
Therefore if $M_n  < n$ and  $N'_n \geq 2 \alpha \tau_n^2$ and
$$
\sum_{k=M_n+1}^n {\bf 1}\{X_k \in \cup_{j=1}^{\lceil 2 \alpha \tau_n^2
\rceil } Q_{n, J'_{n,j}}
\} < \alpha \tau_n^2 ,
$$
then $
\X_{n} \cap Q_{n,J'_{n,j}}
 \neq 
\Po_{n(1-\delta)} \cap Q_{n,J'_{n,j}}
 $ for at most $\lfloor \alpha \tau_n \rfloor$
values of $j \in [1, 2 \alpha \tau_n^2]$, and hence 
$$
N_n \geq
 \sum_{j=1}^{\lceil 2 \alpha \tau_n^2 \rceil}
I_{n,J'_{n,j}} 
 \geq \lceil 2 \alpha \tau_n^2 \rceil -\alpha \tau_n^2 \geq \alpha \tau_n^2.
$$ 
Hence, if  $n(1-2\delta) < M_n < n$
and
  $ N'_n  \geq 2 \alpha \tau_n^2$
and 
$
\sum_{k=M_n+1}^{M_n + \lceil 2 \delta n \rceil }
 {\bf 1}\{X_k \in \cup_{j=1}^{\lceil 2 \alpha \tau_n^2
\rceil } Q_{n, J'_{n,j}}
\} < \alpha \tau_n
$, then $N_n \geq \alpha \tau_n^2$. 
Hence
\bean
P [ N_n <  \alpha \tau_n^2
  | N'_n  \geq 2 \alpha \tau_n^2, n-2 \delta n <
 M_n < n  ]
\\
\leq P [ \Bin (\lceil 2 \delta n \rceil , \lceil 2 \alpha  
\tau_n^2 \rceil  ((\kappa+5) r_n)^d ) >
\alpha  \tau_n^2
].
\eean
Since $nr_n^d$ is assumed bounded, we can choose
$\delta$ small enough so that
the expectation of the binomial variable in the last line is less
than  $(\alpha/2) \tau_n^2$, and then appeal once more
to Lemma \ref{PenL1.1}
 to see that the above conditional probability
decays exponentially in $\tau_n^2$. Combining all these probability
estimates give the desired result. $\qed$ \\

{\em Proof of Theorem \ref{Gthm}.}
Set $p := P[\G(\X_\kappa,1/(\kappa +3)) \sim \Gamma]$.
Let $V_1,V_2,\ldots$ be a sequence of independent Bernoulli
variables with parameter $p$, independent of $\X_n$.
Let
\bean
S'_n :=
\sum_{j=1}^{\min(\lfloor \alpha \tau_n^2 \rfloor,N_n)
}
 {\bf 1} \{\G(\X_n \cap Q_{n,J(n,j)} ; r_n)
 \sim \Gamma\};
~~~~~
Y_n:= G_n-S'_n,
\eean
and
\bean
S_n := S'_n + \sum_{j=1}^{(\lfloor \alpha \tau_n^2 \rfloor - N_n)^+} V_j,
\eean
where $x^+ : = \max(x,0)$ as usual, and the sum $\sum_{j=1}^0$ is
taken to be zero.

For each $j$, given that $I_{n,j}=1$, the distribution
 of the contribution to $G_n$ from points in $Q_{n,j}$ is Bernoulli
with parameter
$P[ \G((\kappa+3)r_n \X_\kappa, r_n) \sim \Gamma]$, which is $p$.
Hence $S_n$ is binomial  $\Bin(\lfloor \alpha
\tau^2_n \rfloor, p).$
Moreover, the conditional distribution
of $S_n,$ given the value of $Y_n$, does not depend on the value
of $Y_n$, and therefore $S_n$ is independent of $Y_n$.
By \eq{normlim1},
\bean
\lfloor \alpha \tau_n^2  \rfloor^{-1/2}
(G_n-\E G_n) \tod \NN (0, \alpha^{-1} \sigma^2 ).
\eean
Moreover,
\bean
\E[|G_n - (Y_n + S_n)| ] =
\E \left[ \sum_{j=1}^{(\lfloor \alpha \tau_n^2 \rfloor - N_n)^+} V_j \right]
\leq p \lfloor \alpha  \tau_n^2 \rfloor P[N_n < \alpha \tau_n^2]
\eean
so that by Lemma \ref{lem1},  both $\tau_n P[G_n \neq Y_n + S_n] $ and
  $\tau_n^{-1}\E[|G_n -  Y_n-S_n|]$ tend to zero as $n \to \infty$.
Hence,
 Theorem \ref{genthm2} (with $\spann_V=1$) is applicable,
with
$ \lfloor \alpha \tau_n^2 \rfloor$
playing the role of $n$ in that result and
$\alpha^{1/2} \tau_n$ playing the role of $ \mmm_n$,
yielding 
\bean
\sup_{k \in \Z }
\left\{ \left| \alpha^{1/2} \tau_n P[G_n =k] -
\alpha^{1/2} \sigma^{-1}
\phi \left(\frac{k - \E G_n}{
(\alpha^{1/2} \tau_n) \alpha^{-1/2} \sigma}
\right)
\right|
\right\}
\to 0,
\eean
as $ n \to \infty$.
Multiplying through by $\alpha^{-1/2}$
 yields \eq{LLT_Ga}.
$\qed$

\section{Proof of Theorem \ref{finranthm}
}
\lbl{secpfSG}
\allco
Recall the definition of $\spann_X$ (the span of $X$) from Section
\ref{secgenresult}.
\begin{lemm}
\lbl{sumlatlem}
 If $X$ and $Y$ are independent
random variables then $\spann_{X+Y}| \spann_X$.
\end{lemm}
{\em Proof.}
If $\spann_{X+Y}=0$ there is nothing to prove. Otherwise, set
 $h= \spann_{X+Y}$.
Then,
considering characteristic functions, observe that
$$
1 = | \E  \exp(2 \pi \imag  (X +Y )/h)|  =
 | \E  \exp(2 \pi \imag  X /h)|  \times
 | \E  \exp(2 \pi \imag  Y/h)|
$$
so that
 $| \E  \exp(2 \pi \imag  X /h)| =1$ and hence $h|\spann_X$. $\qed$ \\

We are in the setup of Section \ref{secstogeo}.
Recall that the point process $\ZZ_n$
consists of $n$ normally distributed marked points 
in $\R^d$,
while $\U_{n,K}$ 
consists of  $n$ uniformly distributed marked points
in $B(K)$.
Set $ \spann_{n,K} := \spann_{H(\U^*_{n,K})} $.
Set $\spann_n:= \spann_{H(\ZZ^*_n)}$,
and recall from \eq{limspandef} that
$\spann(H) := \liminf_{n \to \infty} \spann_n$.

\begin{lemm}
\lbl{lemlat}
Suppose either (i) $H$ has finite range interactions
and $h_{H(\ZZ^*_n)} < \infty$ for some $n$, or (ii)
$H = H^{(\xi)}$ is
induced by a $\kappa$-nearest neighbour functional $\xi(\bx;\X^*)$,
and $h_{H(\ZZ^*_n)} < \infty$ for some $n > \kappa$.
Then
 $\spann(H)  < \infty$, and if
 $\spann(H)  > 0$, 
 there exists $\mu \in \N$
and $K >0$ such that
 $\spann_{\mmu,K} = \spann(H)$.
 If $\spann(H) = 0$, then for any $\eps > 0$
there exists $\mmu \in \N$  and $K >0 $ such that
$\spann_{\mmu,K} < \eps$.
In case (ii), we can take $\mmu$ such that additionally 
$\mmu \geq \kappa +1$.
\end{lemm}
{\em Proof.}
The support of the distribution of $H(\U^*_{n,K})$ is
increasing with $K$,
so
$\spann_{n,K'}|\spann_{n,K}$
for $K' \geq K$.
Hence, there exists a  limit $\spann_{n,\infty} $ such that
\bea
\spann_{n,\infty} =   \lim_{K \to \infty} \spann_{n,K}
\lbl{lanidef}
\eea
and also we have the implication
\bea
\spann_{n,\infty} >0
   \implies
\exists K:
 \spann_{n,K} =
 \spann_{n,\infty}.
\lbl{eqlani}
\eea
Also, for all $K$
the support of the distribution of $H(\U^*_{n,K})$ is
  contained in  the support of $H(\ZZ^*_n)$, so that
\bea
\spann_n  = \spann_{H(\ZZ^*_n)} \leq
 \spann_{n,K}, ~~~
\forall K,
\lbl{1207e}
\eea
  and hence
$
 \spann_{n} \leq \spann_{n,\infty} $
for all $n$.
We assert that in fact
\bea
\spann_{n,\infty} = \spann_n.
\lbl{1207d}
\eea
This is clear
 when $\spann_{n,\infty} =0$.
When $\spann_{n,\infty} > 0$, there exists
a  countable set $S$  with span $\spann_{n,\infty}$ such that
$P[H(\U^*_{n,K}) \in S ] = 1$ for all $K$.
 But then it is easily deduced that $P[H(\ZZ^*_n) \in S] =1$,
so that $\spann_{n} \geq \spann_{n,\infty}$, and combined
with \eq{1207e} this gives \eq{1207d}.

We shall show in both  cases (i) and (ii) that $\spann_n$ tends to a finite
limit;
that is, for both cases we shall show that
\bea
\spann(H) = \lim_{n \to \infty} \spann_n = \lim_{n \to \infty} \spann_{n,\infty}
< \infty.
\lbl{1210a}
\eea
Also, we show in both cases
that
\bea
\spann(H) > 0
\implies \exists n_0 \in \N: \spann_n = \spann(H) ~~ \forall n \geq n_0.
\lbl{1211a}
\eea
If $h(H) >0$, the desired conclusion follows from
\eq{1211a}, \eq{1207d} and \eq{eqlani}.
If $h(H) =0$, the desired conclusion follows
from \eq{1210a} and \eq{lanidef}.

Consider the  case (i), where
 $H$ has finite range interactions. In this case, we
shall show that
 for all  $n$, \bea
\spann_{n+1} | \spann_n, \lbl{ladown}
\eea
and since we assume $\spann_n < \infty$ for some $n$,
\eq{ladown} clearly implies \eq{1210a} and \eq{1211a}.

We now demonstrate \eq{ladown} in case (i) as follows.
By \eq{1207d} and \eq{eqlani},
to prove \eq{ladown} it suffices to prove that
$\spann_{n+1} | \spann_{n,K}$
 for all $K$.
Choose $\tau$ such that  \eq{finraneq} holds.
 There is a strictly positive probability
that the first $n$ points of $\ZZ_n$ lie in $B(K)$ while
the last one lies outside $B(K+\tau)$. Hence by \eq{finraneq}
and translation-invariance,
the support  of the distribution of $H(\ZZ^*_{n+1})$
contains the support of the distribution of $H(\U^*_{n,K}) +
  H(\{(0,T)\})$,  where $T$ is a $\PP_\MM$-distributed
element of $\MM$, independent of $U^*_{n,K}$.  Hence
by Lemma \ref{sumlatlem},
$\spann_{n+1} | \spann_{n,K}$, so
 \eq{ladown} holds as claimed in this case.

Now consider  case (ii), where we assume $H = H^{(\xi)}$ with
$\xi(\bx;\X)$ determined by the $\kappa$ nearest neighbours.
We claim that if
$j \geq \kappa+1$ and $\ell \geq \kappa+1$ then
\bea
\spann_{j+\ell} | \spann_j
~~~{\rm and}~~~
\spann_{j+\ell} | \spann_\ell.
\lbl{1207j}
\eea
By \eq{eqlani} and \eq{1207d}, to
 verify \eq{1207j}
it suffices to show that
\bea
\spann_{j+\ell} | \spann_{j,K} ~~~ \forall K >0.
\lbl{0126b}
\eea
Given $K$, let $B$ and $B'$ be  disjoint balls of radius $K$, distant
 more than $2K$ from each other.
There is a positive probability that $\ZZ_{j+\ell}$ consists
of $j$ points in $B$ and $\ell$ points in $B'$,
and if this happens then (since we assume $\min(j,\ell) >\kappa)$)
 the $\kappa$ nearest neighbours of
the points in $B$ are also in $B$, while
 the $\kappa$ nearest neighbours of
the points in $B'$ are also in $B'$, so that
 $H(\ZZ^*_{j+\ell})$ is the sum of conditionally independent contributions
from the points in
 $B$ and those in $B'$.
Hence the support of the distribution of $H(\ZZ^*_{j+\ell})$ contains
the support of the distribution of $H(\U^*_{j,K})+ H(\tilde{\U}^*_{\ell,K})$,
where $H(\tilde{\U}^*_{\ell,K})$ is defined to be a
 variable with the distribution
of $H(\U^*_{\ell,K})$ independent of $H(\U^*_{j,K})$.
Then \eq{0126b} follows from Lemma \ref{sumlatlem}.

Define
\bean
\spann' = \inf_{n \geq \kappa+1} \spann_n.
\eean
Then
for all $\eps > 0$ we can pick $j \geq \kappa+1$ with
$\spann_j \leq \spann' + \eps$, and then by
\eq{1207j} we have $\spann_\ell \leq \spann' + \eps$
for $\ell \geq j + \kappa +1$. This demonstrates \eq{1210a} for this case
(with $\spann(H) = \spann'$), since we assume $\spann_n < \infty$ for some
$n$.
Moreover, if $\spann(H) > 0,$ then  in the  argument just given we can
take $\eps < \spann(H)$ and then for
 $\ell \geq j + \kappa +1$ we must have
$\spann_\ell |\spann_j$, which can happen
only if $\spann_\ell =\spann_j$,
 so by \eq{1210a},
 in fact $\spann_\ell = \spann_j = \spann(H)$.
That is, we also have \eq{1211a} for this case.
$\qed$ \\

 Since we are in the setting of Section \ref{secstogeo}, we assume (as in
Section \ref{secRGG}) that $f$ is an almost everywhere continuous
probability density function on $\R^d$ with $\fmax < \infty$.
The point process $\X_n \subset \R^d$ is a sample from this density, and
the marked point process $\X_n^* \subset \R^d \times \MM$
 is obtained by giving each point of $\X_n$ a 
$\PP_\MM$-distributed mark. Recall also that
we are given a sequence  $(r_n)$ with $\rho := \lim_{n \to \infty}
n r_n^d \in (0,\infty)$.
 Recall from \eq{0110b} that $H_n(\X^*):= H(r_n^{-1} \X^*)$
for a given translation-invariant $H $.

 Our strategy for proving Theorem \ref{finranthm} goes as follows.
First we choose $\mu,K$ as in Lemma \ref{lemlat}.
Then we choose constants $\beta \geq K$ and $m \geq \mu$
 in a certain way (see below), 
and use the continuity of $f$ to
 pick  $\Theta(n)$ disjoint deterministic balls of
 radius $\beta r_n$ such
that $f$ is positive and almost constant on each of these balls.
 We use a form of rejection sampling to
make the density of points of $\X_n$ in each (unrejected) ball uniform.
We also 
 reject all balls which do not contain exactly
$m$ points of $\X_n$ in a certain
`good' configuration (of non-vanishing probability).
The definition of `good' is chosen in such a way that
 the contribution to $H_n$
from inside an inner ball of radius $K r_n$ is shielded from 
everything outside the outer ball of radius $\beta r_n$.
 We end up with $\Theta(n)$ (in probability) unrejected balls,
and the contributions to $H_n(\X^*_n)$ from the
corresponding inner balls
 are independent (because of the shielding) and
identically distributed (because  of the uniformly distributed
points) so the sum contribution of these inner balls can play the
role of $S_n$ in Theorem \ref{genthm2}.

In the proof of Theorem \ref{finranthm}, we need to
consider certain functions, sets and sequences, defined for $\beta >0$.
For $x \in \R^d$ with $f(x)>0$,
 define the function 
\bea
g_{n,\beta}(x) := \frac{ \inf\{ f(y):y \in B(x;\beta r_n) \} }{
 \sup\{ f(y):y \in B(x;\beta r_n) \} } ,
\lbl{gndef} \
\eea
and for
 $x \in \R^d$ with $f(x)>0$
and
$g_{n,\beta}(x) >0,$ and $  z \in B(x;\beta r_n)$,
define
\bea
p_{n,\beta}(x,z) := \frac{ \inf\{ f(y):y \in B(x;\beta) \} }{
  f(z)  } .
\lbl{pndef}
\eea

Since we assume $f$ is almost everywhere continuous, the function
$g_{n,\beta}$ 
 converges almost everywhere on $\{x:f(x) > 0\}$ to 1.
By Egorov's theorem (see e.g. \cite{Durr}), given $\beta>0$
 there is a set $A_\beta$ with $\int_{A_\beta} f(x) dx \geq 1/2$,
such that $f(x)$ is bounded away from zero on $A_\beta$ and
 $g_{n,\beta}(x) \to 1$ uniformly on $A_\beta$.

Since we assume \eq{rhofin} with $\rho >0$ here, for $n$ large enough $nr_n^d
< 2 \rho$.
Set
$$
\eta(\beta) : =  2^{-(d+2)}  \vold^{-1} \beta^{-d} \fmax^{-1} \rho^{-1}.
$$
Given $\beta>0$, we claim that for  $n$ large enough so that
 $nr_n^d < 2 \rho$, we can (and do) choose  points 
$x_{\beta,n,1},\ldots,x_{\beta,n,\lfloor \eta(\beta) n\rfloor } $
 in $A_\beta$ with
$|x_{\beta,n,j} - x_{\beta,n,k}| > 2 \beta  r_n$ for $1 \leq j<k 
\leq \lfloor \eta(\beta) n\rfloor .$
 To see this we use
a measure-theoretic version of the pigeonhole
principle, as follows.
Suppose inductively that we have chosen $x_{\beta,n,1},\ldots,x_{\beta,n,k}$,
with $k < \lfloor \eta(\beta) n\rfloor .$
Then  let $x_{\beta,n,k+1}$ be the first point, according
to the lexicographic ordering, in the set
$A_{\beta} \setminus \cup_{j=1}^k B(x_{\beta,n,j};2 \beta r_n)$.
This is possible, because  this set is non-empty,
because by subadditivity of measure,
\bean
\int_{\cup_{j=1}^k B(x_{\beta,n,j};2 \beta r_n)}
f(x) dx
\leq k 
\vold (2 \beta r_n)^d \fmax
< \eta(\beta) n \vold (2 \beta r_n)^d \fmax
\\
= nr_n^d/(4 \rho)
< 1/2 \leq \int_{A_{\beta}} f(x) dx,
\eean
justifying the claim.
Define the ball
$$
B_{\beta,n,j} := B(x_{\beta,n,j}, \beta r_n); ~~~~~
B^*_{\beta,n,j} := B(x_{\beta,n,j}, \beta r_n) \times \MM.
$$
 The
balls $B_{\beta,n,1},\ldots,B_{\beta,n,\lfloor \eta (\beta) n\rfloor }$ are disjoint.

Let $W_1,W_2,W_3,\ldots$ be uniformly distributed random variables
in $[0,1]$,
independent of each other and of $(\bX_j)_{j=1}^n$, where
$\bX_j = (X_j,T_j)$. For $k \in \N$, think of
$W_k$ as  an extra  mark attached to the point $X_k$.
 This is used in the rejection sampling procedure.
Given $\beta$, if $X_k  \in B_{\beta,n,j}$,
let us say that the point $X_k$ is
 $\beta$-{\em red} if the associated mark
$W_k$ is less than  $p_{n,\beta}( x_{\beta,n,j}, X_k)$.
 Given that $X_k$ lies in $B_{\beta,n,j}$ and
is $\beta$-red, the conditional distribution
of $X_k$ is uniform over $B_{\beta,n,j}$.

Now let $m \in \N$, and  suppose $\SS$  is a measurable set of configurations
of $m$ points in $B(\beta)$ such that $P[\U_{m,\beta} \in \SS] > 0$.
The number $m$ and the set $\SS$ will be chosen so that
given there are $m$ points of $\X_n$ in ball $B_{\beta,n,j}$, and
given their rescaled configuration of 
lies in the set $\SS$, there is a subset of these
$m$ points which are `shielded' from the rest of $\X_n$.

Given $\SS$ (and by implication $\beta$ and $m$),
for $1 \leq j \leq \lfloor \eta(\beta) n\rfloor $,
let $I_{\SS,n,j}$ be the indicator of the event that
the following conditions  hold:
\begin{itemize}
\item
The point set
  $\X_n \cap B_{\beta,n,j}$ consists of $m$ points, all of them
$\beta$-red;
\item The configuration $r_n^{-1}(-x_{\beta,n,j} +
 ( \X_n \cap B_{\beta,n,j} ))  $ is in $\SS$.
\end{itemize}
Let $N_{\SS,n} : = \sum_{j=1}^{\lfloor \eta (\beta) n\rfloor } I_{\SS,n,j}$,
 and list the $i$ for which $I_{\SS,n,j} =1$ in increasing order
as $J(\SS,n,1)  \ldots,$ $ J(\SS,n,N_{\SS,n})$.

\begin{lemm}
\lbl{redlem}
Let $\beta > 0 $,  and $m \in \N$.  Let $\SS$  be a measurable set of
 configurations of $m$ points in $B(\beta)$ such that
 $P[\U_{m,\beta} \in \SS] > 0$.
Then: (i) there  exists $\delta > 0$ such that
\bea
\limsup_{n \to \infty}
\left(
n^{-1} \log
P [N_{\SS,n} < \delta n] \right) < 0,
\lbl{1207c}
\eea
and (ii) conditional on
the values  of $I_{\SS,n,i}$ for $1 \leq i \leq \lfloor \eta(\beta)n\rfloor$
and the configuration of $\X_n$
outside $B_{\beta,n,J(\SS,n,1)} \cup \cdots \cup
B_{\beta,n,J(\SS,n,N_{\SS,n})}$,
the joint distribution of the point sets
$$
r_n^{-1}(-x_{\beta,n,J(\SS,n,1)} +( \X_n \cap B_{\beta,n,J(\SS,n,1)}) )
, \ldots,
r_n^{-1} (-x_{\beta,n,J(\SS,n,N_{\SS,n})} +
  (\X_n \cap B_{\beta,n,J(\SS,n,N_{\SS,n})}
))
$$
 is that of
$N_{\SS,n}$
 independent copies of $\U_{m,\beta}$ each conditioned
to be in $\SS $.
\end{lemm}
{\em Proof.}
Consider first the asymptotics for $\E [N_{\SS,n}]$.
Given a finite point set $\X \subset \R^d$
 and a set $B \subset \R^d$, let $\X(B)$ denote
the number of points of $\X$ in $B$.
Fix $m$.
Since
$f$ is bounded away from zero and infinity on $A_\beta$ and $g_{n,\beta} \to 1$
uniformly on $A_\beta$,
we have
uniformly over   $x \in A_\beta$ that 
$$
n \int_{B(x; \beta r_n)} f(y)dy = n f(x)
 \int_{B(x; \beta r_n)} (f(y)/f(x))dy \to \beta^d \vold \rho  f(x)
$$
Hence by
 binomial approximation to Poisson,
\bean
P [\X_n(B(x;\beta  r_n)) = m ] \to
 \frac{ ( \beta^d \vold  \rho f(x) )^m
 \exp( -  \beta^d \vold \rho f(x) ) }{ m! }
~~~ {\rm as} ~n \to \infty,
\eean
and this convergence is also uniform over $x \in A_\beta$.

Given $m$ points $X_k$ in $B_{\beta,n,j}$, the probability
that these are all $\beta$-red is  at least
$g_{n,\beta}(x)^m$ so  exceeds
$\frac{1}{2}$
if $n$ is large enough, since
 $g_{n,\beta} \to 1$
uniformly on $A_\beta$.

Given that $m$ of the points $X_k $ lie in $ B_{\beta,n,j} $,
and given that they are all  $\beta$-red, their spatial locations are
independently uniformly distributed over
 $ B_{\beta,n,j}$; hence the conditional probability that
$r_n^{-1}(- x_{\beta,n,j}  +  ( \X_n \cap B_{\beta,n,j} ) )$ lies in $ \SS$
is a strictly positive constant.

These arguments show that
  $ \liminf_{n \to \infty} n^{-1} \E [ N_{\SS,n} ] > 0$.
They also demonstrate part (ii) in the statement of
the lemma.

Take $\delta > 0$ with 
$2 \delta < \liminf_{n \to \infty} n^{-1} \E [N_{\SS,n} ] $.
We shall show that $P[N_{\SS,n} < \delta n]$ decays exponentially in $n$,
using Lemma \ref{Azulem}.
The variable $N_{\SS,n}$ is a function of $n$ independent
identically distributed
triples  (marked points) $(X_k,T_k,W_k)$.

Consider the effect of changing the value of one of the marked
points ($(X,T,W)$ to $(X',T',W')$, say). The change could
affect the value of $I_{\SS,n,j}$
for at most two values of $j$,
namely the $j$  with $X \in B_{\beta,n,j}$ and the $j'$  with
$X' \in B_{\beta,n,j'}$.
So by Lemma \ref{Azulem}, 
$$
P[ | N_{\SS,n} - \E N_{\SS,n} | > \delta n]
\leq 2 \exp ( - \delta^2 n / 8),
$$
 and \eq{1207c} follows.
$\qed$ \\

{\em Proof of Theorem \ref{finranthm} under condition (i) (finite range
interactions).}  Recall that $\spann(H) $ is given by
\eq{limspandef}. 
Since condition (i) includes the assumption
that $\spann_{H(\ZZ^*_n)} < \infty$ for some
$n$, by Lemma \ref{lemlat}
we have $\spann(H)< \infty$.
Let $b>0$ with $\spann(H)|b$. Let $\eps \in ( 0,b)$.
Let $\mmu \in \N$, and $K >0$, be 
as given by Lemma \ref{lemlat}.
Then  $\spann_{\mmu,K} = \spann(H)$ if $\spann > 0$, or
 $\spann_{\mmu,K} < \eps $ if $\spann = 0$.
 Moreover
$ H(\U^*_{\mmu,K})$ is integrable by
assumption.
Set
\bea
b_1 :=
  \begin{cases}
\spann_{\mmu,K} \lfloor b/\spann_{\mmu,K} \rfloor &
 \textrm{~if~} \spann_{\mmu,K} > 0 \\
b
 & \textrm{~if~} \spann_{\mmu,K} = 0 .
\end{cases}
\lbl{b1def}
\eea

Choose $\tau \in (0, \infty)$ such that \eq{finraneq} holds.
We shall apply Lemma \ref{redlem} with
 $\beta = K+ \tau$. 
 Let $\SS$ be the set of configurations
of $\mmu$ points in $B(K+\tau)$ such that
in fact all of the points are in $B(K)$.
By Lemma \ref{redlem},
we can find $\delta > 0$ such that, writing $N_n$ for  $N_{\SS,n}$
we have exponential decay of $P[N_n < \delta n]$.

Let $V_1,V_2,\ldots,$ be random variables distributed as independent
 copies of $H(\U^*_{\mmu,K})$, independently of $\X^*_n$.
Set
\bean
S'_n :=
\sum_{\ell=1}^{\min(\lfloor \delta n \rfloor,N_n) } H_n(
\X^*_n \cap B^*_{K+\tau,n,J(\SS,n,\ell)}  ); ~~~~
S_n = S'_n + \sum_{j=1 }^{(\lfloor \delta n \rfloor -N_n)^+} V_j.
\eean
Thus, $S'_n$ is the the total contribution to $H_n(\X^*_n)$ from points
in $\cup_{\ell=1}^{ \min(\lfloor \delta n\rfloor
 ,N_n) } B^*_{K+\tau,n,J(\SS,n,\ell)}$.

By Part (ii) of Lemma \ref{redlem},
given that $N_n \geq \delta n$, for each $\ell$
  we know that $r_n^{-1}(- x_{\beta,n,J(\SS,n,\ell)} + \X^*_n ) \cap
 B^*(K+\tau)$
is conditionally distributed as $\U^*_{\mmu,K+ \tau}$ conditional
on $\U^*_{\mmu,K+\tau} \in \SS$; in other words, distributed
as $\U^*_{\mmu,K}$.  Therefore the distribution of $S_n$
is that of the sum of  $\lfloor \delta n\rfloor$
 independent copies of $H(\U^*_{\mmu,K})$,
independent of the contribution of the other points.
Let $Y_n$ denote the contribution of the other points, i.e.
\bean
Y_n:= H_n (\X_n^* ) -  S'_n.
\eean
Since the distribution
of $S_n,$ given the value of $Y_n$, does not depend on the value
of $Y_n$, $S_n$ is independent of $Y_n$.

By assumption $H_n(\X^*_n)$ and $S_n$ are
integrable.
Clearly $n^{1/2}P[H_n(\X_n) \neq Y_n + S_n]$ is at most
$n^{1/2}P[N_{n} < \delta n]$,
which tends to zero by \eq{1207c}.
Also by conditioning on $N_{n}$, we have that
\bea
n^{-1/2}\E [|H_n(\X^*_n) - (Y_n + S_n) |]
= n^{-1/2}\E  \left[ \left|\sum_{j=1}^{(\lfloor \delta n \rfloor
 - N_{n})^+} V_j \right| \right]
\nonumber \\
\leq n^{-1/2} \E [ (\lfloor \delta n \rfloor - N_n)^+ ]
\E \left[ \left| V_1 \right| \right]
\nonumber \\
\leq n^{-1/2} \lfloor \delta n \rfloor
P[ N_{n} \leq \delta n] \E \left[ \left| V_1 \right| \right],
\label{0113e}
\eea
which tends to zero by \eq{1207c}.
This also shows that $Y_n$ is integrable
By the assumption \eq{Hclteq},
\bea
 \lfloor\delta n\rfloor^{-1/2}
(H_n(\X^*_n)-\E H_n(\X^*_n)) \tod \NN (0,  \delta^{-1} \sigma^2 ),
\lbl{lim3}
\eea
and so, since $\spann_{\mmu,K}|b_1$,
Theorem \ref{genthm2} is applicable, and
  yields
\bea
\sup_{u \in \R }
\left\{ \left| (\delta n)^{1/2} P[H_n(\X^*_n) \in [u,u+b_1)] -
\delta^{1/2} \sigma^{-1} b_1 \phi  \left( \frac{u- \E H_n(\X^*_n) }{(\delta
n)^{1/2} ( \delta^{-1} \sigma^2)^{1/2}}
  \right)
\right|  \right\} \to 0,
\nonumber \\
\label{0113d}
\eea
and dividing through by $\delta^{1/2}$ gives \eq{1205a}
in all cases where $b = b_1$. 
 In general, suppose $b \neq b_1$.
Then $\spann(H) =0$ (else $\spann_{\mmu,K}= \spann(H)$ and $\spann(H) |b$
so $b = b_1$ by \eq{b1def}), and hence
$\spann_{\mmu,K} < \eps$.
 Since $b_1 \leq b$ by \eq{b1def},
we have
that
\bean
\inf_{u \in \R }
\left\{
  n^{1/2} P[H_n(\X^*_n) \in [u,u+b)] -
 \sigma^{-1} b \phi  \left( \frac{u- \E H_n(\X^*_n) }{
n^{1/2} \sigma}
  \right)
\right\}
\\
\geq
\inf_{u \in \R }
\left\{
  n^{1/2} P[H_n(\X^*_n) \in [u,u+b_1)] -
 \sigma^{-1} b_1
\phi  \left( \frac{u- \E H_n(\X^*_n) }{
n^{1/2} \sigma}
  \right)
\right\}
\\
+ \sigma^{-1}(b_1-b)
(2 \pi)^{-1/2}
\eean
so that by \eq{0113d}, since $b_1 \geq b-\eps$,
\bean
\liminf_{n \to \infty}
\inf_{u \in \R }
\left\{
  n^{1/2} P[H_n(\X^*_n) \in [u,u+b)] -
 \sigma^{-1} b \phi  \left( \frac{u- \E H_n(\X^*_n) }{
n^{1/2} \sigma}
  \right)
\right\}
\geq - \frac{\eps}{\sigma} (2 \pi)^{-1/2}.
\eean
Similarly, setting
 $b_2 : = \spann_{\mmu,K} \lceil b/\spann_{\mmu,K} \rceil$,
we have that
\bean
\sup_{u \in \R }
\left\{
  n^{1/2} P[H_n(\X^*_n) \in [u,u+b)] -
 \sigma^{-1} b \phi  \left( \frac{u- \E H_n(\X^*_n) }{
n^{1/2} \sigma}
  \right)
\right\}
\\
\leq
\inf_{u \in \R }
\left\{
  n^{1/2} P[H_n(\X^*_n) \in [u,u+b_2)] +
 \sigma^{-1} b_2
\phi  \left( \frac{u- \E H_n(\X^*_n) }{
n^{1/2} \sigma}
  \right)
\right\}
\\
+ \sigma^{-1}(b_2-b)
(2 \pi)^{-1/2}
\eean
so that since $b_2 -b \leq \eps$,
\bean
\limsup_{n \to \infty}
\sup_{u \in \R }
\left\{
  n^{1/2} P[H_n(\X^*_n) \in [u,u+b)] -
 \sigma^{-1} b \phi  \left( \frac{u- \E H_n(\X^*_n) }{
n^{1/2} \sigma}
  \right)
\right\}
\leq  \frac{\eps}{\sigma} (2 \pi)^{-1/2}.
\eean
Since $\eps >0 $ is arbitrarily small, this gives us \eq{1205a}.
$\qed$ \\

{\em Proof of Theorem \ref{finranthm} under condition (ii).}
We now assume
that $H$, instead of having finite range,
 is given by \eq{induceh} with $\xi$ depending
only on the $\kappa$ nearest neighbours.
Again, by Lemma \ref{lemlat} we have that $\spann(H) $,
given by \eq{limspandef}, is finite.

Let $b>0$ with $\spann(H)|b$. Let $\eps \in ( 0,b)$.
Let $\mmu  \in \N$ and $K >0$, with $\mu \geq \kappa +1$,
by   as given by
Lemma \ref{lemlat}.
Then  $\spann_{\mmu,K} = \spann(H)$ if $\spann(H) > 0$, and
 $\spann_{\mmu,K} < \eps $ if $\spann(H) = 0$.
Also, $ H(\U^*_{\mmu,K}) $
 integrable, by the integrability assumption
in the statement of the result being proved.

Let $\BB_1,\BB_2,\ldots, \BB_{\nu}$ be a minimal collection of
open balls of radius $K$, each of them centred at a point on
  the boundary of $B(4K)$, such that their union contains
  the boundary of $B(4K)$.  
 Let $\BB_0$ be the ball $B(K)$.

We shall apply  Lemma \ref{redlem} with $\beta= 5K$,
with $m = (\nu +1) \mmu$, and with $\SS$ as follows.
$\SS$ is the set of configurations of $m =(\nu +1) \mmu$
points in $B(\beta) = B(5K)$, such that each of
$\BB_1,\ldots,\BB_{\nu}$ contains  at least $\mmu$ points,
and $\cup_{i=1}^\nu \BB_i$ contains exactly $\nu \mmu$ points,
 and also the ball $\BB_0$ contains exactly $\mmu$ points
(so that consequently there are no points in $B(5K)
\setminus \cup_{i=0}^\nu \BB_i$).
A similar construction (using squares rather than balls, and with diagram)
was given by Avram and Bertsimas \cite{AB} for a related problem.

With this choice of $\beta$ and $\SS$, let the locations
$x_{\beta,n,j} = x_{5K,n,j}$, the balls $B_{\beta,n,j} = B_{5K,n,j}$,
 the indicators
 $I_{\SS,n,j}$, and the variables $N_{\SS,n}$ and
$J(\SS,n,\ell)$  be as described just before
 Lemma \ref{redlem}.
By that result, we can (and do) choose $\delta > 0$ such that
\eq{1207c} holds.

For $ 1 \leq \ell \leq N_{\SS,n}$, the point process
 $r_n^{-1}( - x_{5K,n,J(\SS,n,\ell)} + (\X_n \cap B_{5K,n,J(\SS,n,\ell)})) $
  has $\mmu $ points within distance
$K$ of the origin, and
also at least $\mmu$ points in each of the balls $\BB_1,\ldots,\BB_\nu$.

Since $\mmu \geq \kappa +1$, for any point configuration in $\SS$,
each point inside $B(K)$ has its $\kappa$ nearest neighbours
also inside $B(K)$.  Also none of the  points in $B(5K) \setminus B(K)$
has any of its $\kappa$ nearest neighbours in $B(K)$.
Finally, any further added point outside $B(5K)$ cannot have
any of its $\kappa$ nearest neighbours inside $B(K)$, since
the line segment from such a point to any point in
 $B(K)$ passes through the boundary of $B(4K)$ at a location
inside some
$\BB_i$, and any of the $\mmu$ or more points inside $\BB_i$ are
closer to the outside point than the point in $B(K)$ is.
To summarise this discussion, the points in $B(K)$ are shielded from those
outside $B(5K)$.

Given $n$, let
 $\WW_{(\nu+1)\mmu,5K}^{(1)},\ldots,
\WW_{(\nu+1)\mmu,5K}^{(\lfloor \delta n \rfloor )}$
 be a collection of (marked) point processes which are
each distributed as $\U^*_{(\nu+1)\mmu,5K}$ conditioned on 
$\U^*_{(\nu+1)\mmu,5K} \in \SS$,
 independently of each other and of $\X^*_n$.
For  $1 \leq j \leq \lfloor \delta n\rfloor $
set $V_j := H(\WW_{(\nu+1)\mmu,5K}^{(j)} \cap B^*(K))$,
so that $V_1,V_2,\ldots$ $V_{\lfloor \delta n\rfloor }$ are
  random variables distributed as
  independent copies
 of $H(\U^*_{\mmu,K})$, independent of $\X_n$.
Define $S'_n$ and $S_n$ by
\bean
S'_n :=
\sum_{\ell=1}^{\min(\lfloor \delta n \rfloor,N_{\SS,n}) }
H_n( \X^*_n \cap B^*(x_{5K,n,J(\SS,n,\ell)},K r_n)  ) ; ~~~
S_n := S'_n + \sum_{j=1}^{( \lfloor \delta n \rfloor - N_{\SS,n})^+ } V_j.
\eean
Also set
$
Y_n:= H_n(\X^*_n) -S'_n.
$

Thus $S'_n$ is the total contribution to $H_n(\X^*_n)$
from points in $B^*(x_{5K,n,J(\SS,n,\ell)};Kr_n)$,
$1 \leq \ell \leq \min(\lfloor \delta n\rfloor ,N_{\SS,n})$.
On account of the shielding effect described above,
 $S_n$ is the
sum of $\lfloor \delta n \rfloor $ independent copies of a random variable
with the distribution of $H(\U^*_{\mmu,K})$.
Moreover, we assert  that the distribution
of $S_n,$ given the value of $Y_n$, does not depend on the value
of $Y_n$, and therefore $S_n$ is independent of $Y_n$.

Essentially, this assertion holds  because
for any triple of sub-$\sigma$-algebras
 $\F_1,\F_2,\F_3$, if $\F_1 \vee \F_2$ is independent of
$\F_3$ and $\F_1$ is independent of $\F_2$ then $\F_1$ is independent of
$\F_2 \vee \F_3$ (here $\F_i \vee \F_j$ is the smallest $\sigma$-algebra
containing both $\F_i$ and $\F_j$).
  In the present instance, to define these $\sigma$-algebras we
first define the marked point processes $\Y_j $ for
$1 \leq j \leq \lfloor \delta n\rfloor$
by
$$
\Y_j :=
  \begin{cases}
r_n^{-1} ( -x_{5K,n,J(\SS,n,j)} + (\X^*_n  \cap B^*_{5K,n,J(\SS,n,j)} ))
 &
 \textrm{~if~} 1 \leq j \leq \min (\lfloor \delta n \rfloor , N_{\SS,n} )
 \\
\WW^{(j- N_{\SS,n})}_{(\nu +1)\mmu,5K}
 & \textrm{~if~}
N_{\SS,n} < j \leq \lfloor \delta n \rfloor .
\end{cases}
$$
Take $\F_3$ to be the
$\sigma$-algebra generated by
the values of $J(\SS,n,1), \ldots,$
\linebreak
$ J(\SS,n,\min(\lfloor \delta n\rfloor ,N_{\SS,n}))$ and  the
 locations and marks of  points of $\X_n$ outside the union of the balls
$B_{5K,n,J(\SS,n,1)}, \ldots,$ $ B_{5K,n,J(\SS,n,\min(\lfloor
\delta n \rfloor,N_{\SS,n}))}$.
Take $\F_2$ to be the
$\sigma$-algebra generated by the point processes
$ \Y_{j} \cap B^*(5K) \setminus B^*(K), 1 \leq j \leq \lfloor \delta n \rfloor$.
Take $\F_1$ to be the
$\sigma$-algebra generated by the point processes
$\Y_{j} \cap B^*(K), 1 \leq j \leq \lfloor \delta n \rfloor$.
Then by Lemma \ref{redlem} and the definition of $\SS$,
 $\F_1 \vee \F_2$ is independent of
$\F_3$ and $\F_1$ is independent of $\F_2$, so $\F_1$ is independent of
$\F_2 \vee \F_3$.
The variable $S_n$ is measurable with respect to $\F_1$, and
by shielding, the variable $Y_n$ is measurable with respect to $\F_2 \vee \F_3$,
justifying our assertion of independence.

By the assumptions of the result being proved,
 $H_n(\X^*_n)$ and $S_n$ are integrable.
Clearly $n^{1/2}P[H_n(\X^*_n) \neq Y_n + S_n]$ is at most
$n^{1/2}P[N_{\SS,n} < \delta n]$,
which tends to zero. Also, as
with \eq{0113e} in Case (i), we have that
$n^{-1/2}\E [|H_n(\X^*_n) - (Y_n + S_n) |]$
 tends to zero by \eq{1207c}, and  $Y_n$ is
integrable.  By \eq{Hclteq},
\bea
 \lfloor \delta n\rfloor^{-1/2}
(H_n(\X^*_n) -\E H_n(\X^*_n) ) \tod \NN (0,  \delta^{-1} \sigma^2 ),
\eea
and so, since  $\spann_{\mmu,K}|b_1$,
 Theorem \ref{genthm2} is applicable
with $Z_n = H_n(\X^*_n)$,  yielding
\bean
\sup_{\{u \in \R \}}
\left\{ \left| (\delta n)^{1/2} P[H_n(\X^*_n) \in [u,u+b_1)] -
\delta^{1/2} \sigma^{-1} b_1
\phi \left(\frac{u - \E H_n(\X^*_n)}{
(\delta n)^{1/2} \delta^{-1/2} \sigma}
\right)
\right|
\right\}
\to 0,
\eean
as $ n \to \infty$.
Multiplying through by $\delta^{-1/2}$
yields \eq{1205a} for this case, when $b_1 =b$. If $b_1 \neq b$,
we can complete the proof in the same manner as in the proof
for Case (i). $\qed$

\section{Proof of Theorems 
 \ref{fin0theo}, \ref{fintheo} and \ref{nnlclt}}
\lbl{secusePEJP}

\allco

The proofs of Theorems \ref{fin0theo}, \ref{fintheo} and \ref{nnlclt}
 all rely heavily on
Theorem 2.3 of \cite{PenEJP} so for convenience we state
that result here in the form we shall use it.
This requires some further notation, besides the notation
we set up earlier in Section \ref{secstogeo}.

As before, we assume $\xi(\bx,\X^*)$ 
is a
 translation invariant,
  measurable $\R$-valued function
defined for all pairs
 $(\bx,\X^*)$, where $\X^* \subset \R^d\times \MM$
 is finite  and $\bx$
is an element of $\X^*$. 
We extend the definition of $\xi(\bx,\X^*)$ 
to the case where $\X^* \subset \R^d\times \MM$
and $\bx \in ( \R^d \times \MM) \setminus \X^*$,
by setting $\xi(\bx,\X^*)$ 
to be $\xi(\bx,\X^* \cup \{\bx\})$ in this case.
Recall that
$H^{(\xi)}$ is defined by \eq{induceh}.

Let $T$ be an $\MM$-valued random variable with 
distribution $\PP_\MM$, independent of everything else.
  For $\lla >0$ let
  $M_{\lla}$ be a Poisson variable with parameter $\lla$,
independent of everything else, and  let $\Po_{\lla}$ be the
 point process $\{X_1,\ldots,X_{M_{\lla}}\}$, which is a Poisson
point process with intensity $\lla f(\cdot)$. 
Let $\Po_\lla^* := \{(X_1,T_1),\ldots,
(X_{M_{\lla}},T_{M_\lla})\}$ be the corresponding
marked Poisson process.

Given $\lla >0$,
we say $\xi$ is $\lla$-homoegeneously stabilizing if
there is an almost surely finite positive random variable $R$ such that
with probability 1, 
$$
\xi((0,T); (\H_\lla^* \cap B^*(0;R)) \cup {\cal Y} )
=
\xi((0,T); \H_\lla^* \cap B^*(0;R)) 
$$
for all finite ${\cal Y} \subset (\R^d \setminus  B(0;R)) \times \MM$.
Recall that $\supp(f)$ denotes the support of $f$.
We say that $\xi $ is {\em exponentially stabilizing} if
for $\lla \geq 1$ and $x \in \supp(f)$  there exists
a  random variable $R_{x,\lla}$ such that
\bean
\xi((\lla^{1/d}x,T);
\lla^{1/d} (\Po_{\lla}^* 
\cap B^*(x;\lla^{-1/d} R_{x,\lla} )) \cup {\cal Y} )
\\
= 
\xi((\lla^{1/d}x,T); 
\lla^{1/d} ( \Po_{\lla}^* \cap B^*(x;\lla^{-1/d} R_{x,\lla} ) ))
\eean
for all finite ${\cal Y} \subset (\R^d \setminus  B(x;
\lla^{-1/d} R_{x,\lla}) )) \times \MM$,
and there exists a finite positive
constant  $C $  such that
$$
P[ R_{x,\lla} >s] \leq C \exp(-C^{-1}s), ~~~ s \geq 1, ~
\lla \geq 1,~ f \in \supp(f).
$$
For $k \in \N \cup \{0\} $,
let $\TT_k$ be the collection of all
subsets of $\supp(f)$ with at most $k$ elements.
For $k \geq 1$ and
 $\A =\{x_1,\ldots,x_k\} \in  \TT_k \setminus \TT_{k-1}$,
let $\A^*$ be the corresponding marked point set
$\{(x_1,T_1),\ldots,(x_k,T_k)\}$ where
$T_1,\ldots,T_k$ are independent $\MM$-valued
variables with distribution $\PP_\MM$, independent of everything else.
If $\A \in \TT_0$ (so $\A = \emptyset$)
let $\A^*$ also be the empty set. 

We say that $\xi$ is {\em binomially exponentially
stabilizing} if
there exist finite positive constants  $C,\eps $ such that 
 for all $x \in \supp(f)$ and
all $\lla \geq 1 $ and
 $n \in \N \cap ( (1-\eps) \lla, (1+\eps)\lla)) $,
and $\A \in \TT_2$,
 there is a random variable $R_{x,\lla,n,\A}$  such that
\bea
\xi((\lla^{1/d}x,T);
\lla^{1/d} ((\X_n^* \cup \A^*) 
\cap B^*(x;\lla^{-1/d} R_{x,\lla,n,\A} )) \cup {\cal Y} )
\nonumber \\
= 
\xi((\lla^{1/d}x,T); 
\lla^{1/d} ( (\X_{n}^* \cup \A^*) \cap B^*(x;\lla^{-1/d} R_{x,\lla,n,\A} ) ))
\lbl{BiRS}
\eea
for all finite ${\cal Y} \subset (\R^d \setminus  B(x;
\lla^{-1/d} R_{x,\lla,n,\A}) ) \times \MM$, and such that
all $\lla \geq 1$  
and all $n \in \N \cap ( (1-\eps) \lla, (1+\eps)\lla)) $,
and all $x \in \supp(f)$ and all $\A \in \TT_2$,
$$
P[ R_{x,\lla,n,\A} >s] \leq C \exp(-C^{-1}s), ~~~ s \geq 1.
$$
Given 
$p >0$ and
$\eps >0$, we 
 consider the moments conditions
\bea
\lbl{Pomoments}
\sup_{\lla \geq 1, x \in \supp(f), \A \in \SS_1} \E[
| \xi((\lla^{1/d}x,T); \lla^{1/d} (\Po^*_\lla  \cup \A^*  ))|^p]
 < \infty
\eea
and 
\bea
\label{Bimoments}
\sup_{\lla \geq 1, n \in \N \cap ((1-\eps)\lla,(1+\eps)\lla) ,
x \in \supp(f), \A \in \SS_3} \E[|
 \xi((\lla^{1/d}x,T); \lla^{1/d} ( \X_n^* \cup \A^*  ))|^p]
 < \infty.
\eea
\begin{theo}
\lbl{lemPEJP}
Suppose $H = H^{(\xi)} $ is induced by translation-invariant
 $\xi$. Suppose that $\xi$ is
$f(x)$-homogeneously 
stabilizing for Lebesgue-almost all $x \in \supp(f)$, and
$\xi$ is exponentially stabilizing, binomially exponentially
stabilizing and for some $\eps >0$ and $p >2$
satisfies \eq{Pomoments}
and \eq{Bimoments}. 
Suppose $\fmax < \infty$ and $\supp(f)$ is bounded.
Suppose $(\lla(n),n \geq 1)$
is a sequence taking values in $\R^+$
with
 $|\lla(n)-n| = O(n^{1/2})$
as $n \to \infty$. Then  
there exists $\sigma \geq 0$ such that
$$
n^{-1/2} (H^{(\xi)}( \lla(n)^{1/d} \X^*_n)
- \E H^{(\xi)} (\lla(n)^{1/d} \X^*_n)) \tod \NN (0,\sigma^2),
$$
and 
$n^{-1}{\rm Var}(H^{(\xi)}(\lla(n)^{1/d} \X^*_n) \to \sigma^2$ 
as $n \to \infty$.

\end{theo}
Theorem \ref{lemPEJP} is a special case of Theorem 2.3 of \cite{PenEJP},
which also provides an expression for $\sigma$
in terms of integrated two-point correlations; 
that paper considers random measures
given by a sum of contributions from each point,
whereas here we just consider the total measure.  The sets 
 $\Omega_\infty$ and (for all $\lla \geq 1$) $\Omega_\lla $ 
in \cite{PenEJP} are taken to be 
$\supp(f)$. Our $\xi$ is translation invariant,
and these assumptions lead to some simplification
of the notation in \cite{PenEJP}.  \\

{\bf Proof of Theorem \ref{fin0theo}.}
The condition that $\xi(\bx;\X^*)$ has finite range
implies that $H = H^{(\xi)} $ has finite range interactions.
Since $\xi$ has finite range $r$,
 $\xi$ is $\lla$-homogeneously stabilizing for all $\lla >0$, 
exponentially stabilizing and binomially exponentially
stabilizing 
(just take $R = r$, $R_{x,\lla}=r$ and $R_{x,\lla,n,\A}=r$).

We shall establish \eq{Hclteq} by applying 
Theorem \ref{lemPEJP}.  We need to check the
   moments conditions \eq{Pomoments} and \eq{Bimoments}
 in the present setting.
Since we assume that $\fmax< \infty$, 
for any $\lambda >0$ and  any  $n \in \N$ with
 $n \leq 2 \lla$, and any $x \in \supp(f)$, the   variable
 $\card(\X_n^* \cap B^*(x;r \lla^{-1/d}))$ 
is binomially distributed with with mean  at most $\omega_d \fmax 2 r^d$.
Hence by Lemma \ref{PenL1.1}, there is a 
constant $C$, such that whenever $n \leq 2 \lla$ and
$x \in \supp(f)$ we have
\bea
P[ \card (\X_n^* \cap B^*(x;r \lla^{-1/d})) > u ] 
\leq  C \exp(- u/C), ~~~ u \geq 1.
\lbl{0721a}
\eea 
Moreover by \eq{polyxi} 
and the assumption that $\xi$ has range $r$,
for $\A \in \TT_3$
we have
\bean
\E [\xi((\lla^{1/d}x,T);\lla^{1/d}(\X^*_n \cup {\cal A}^* ) )^4]
\leq \polybeta^4 \E[ (4+ \card (\X_n^* \cap B^*(x;r \lla^{-1/d})))^{4 \polybeta}  ] 
\eean
so by \eq{0721a} we can bound the fourth moments of
$\xi((\lla^{1/d}x,T);\lla^{1/d}(\X^*_n \cup {\cal A}^* ) )$
uniformly over  $(x,\lla,n, \A)  \in \supp(f) \times [1,\infty)
\times \N \times  \TT_3$ with $n \leq 2 \lla$.
This gives us \eq{Bimoments} (for $p=4$ and $\eps = 1/2$)
 and \eq{Pomoments}
 may be  deduced similarly.

Hence, the assumptions 
 of Theorem \ref{lemPEJP} are satisfied,  with $\lambda(n)$
in that result  given by $\lambda(n) = r_n^{-d}$.
By Theorem \ref{lemPEJP}, for some
$\sigma \geq 0$ we have \eq{Hclteq} and
\eq{varconv}. Then by Theorem \ref{finranthm}, 
 we can deduce that $\sigma >0$ and
 $\spann(H) < \infty$  and 
 \eq{1205a} holds whenever $\spann(H) | b$.
$\qed$ \\

{\bf Proof of Theorem \ref{fintheo}.}
%
 Under condition \eq{finraneq}, the
functional $H(\X^*)$ can be expressed as
a sum of contributions from components of the geometric
(Gilbert) graph  $\G(\X,\tau)$,  where $\X := \pi(\X^*)$
 is the unmarked point set
corresponding to $\X^*$
 (recall that $\pi$ denotes the canonical projection from 
$\R^d \times \MM$ onto $\R^d$.)
Hence, $H(\X^*)$ can be written as $H^{(\xi)}(\X^*)$
where $\xi(\bx;\X^*)$ denotes the contribution to $H(\X^*)$ from the
component containing $\pi(\bx)$, divided by the number of vertices
in that component.
Then $\xi(\bx;\X^*)$ is unaffected by changes to $\X^*$ that do not
affect the component of $\G(\X,\tau)$ containing $\pi(\bx)$, and
we shall use this to demonstrate that the
conditions of Theorem \ref{lemPEJP} hold,
as follows
(the argument is similar to that in Section 11.1 of \cite{Penbk}).

Consider first the homogeneous stabilization
 condition.
For $\lla >0$, let $R(\lla)$ be 
the  maximum Euclidean distance from
the origin of vertices in the graph
$\G(\H_{\lla}\cup \{0\}, \tau)$
that are pathwise connected to the origin,
which by scaling (see the Mapping theorem in
\cite{Kingman})  has the same distribution as
$\tau $ times the  maximum Euclidean distance from
the origin of vertices in $\G(\H_{\tau^d \lla} \cup \{0\},1)$,
that are pathwise connected to the origin.
Then 
  $R(\lla)$ is almost surely finite, for any  
$\lla \in (0,\tau^{-d}\lla_c)$.

Changes to $\H_{\lla}$ at a distance more than $R(\lla) + \tau$
from the origin do not affect the component of
$\G(\H_{\lla} \cup \{0\},\tau)$ containing
the origin and therefore do not affect $\xi((0,T);\H^*_{\lla})$.
This shows that 
 $\xi$ is $\lla$-homogeneously stabilizating
for any $\lla < \tau^{-d} \lla_c$, and therefore
by assumption \eq{subcrit}
the homogeneous stabilization condition of Theorem \ref{lemPEJP} 
holds.

Next we consider the binomial stabilization condition.
 Let $x \in \supp(f) $.
Let $R_{x,\lambda,n}$ be equal to $\tau$ plus
the maximum Euclidean  distance from $\lambda^{1/d}x$
of vertices in $\G(\lambda^{1/d}(\X_n \cup \{x\}),\tau)$ 
that are pathwise connected to $\lambda^{1/d}x$.
Changes to $\X_n$ at a Euclidean distance greater than
$\lambda^{-1/d} R_{x,\lambda,n}$ 
from $x$  will have no effect on $\xi((\lambda^{1/d}x,T);
\lambda^{1/d}\X_n^* )$.

Using \eq{subcrit}, let $\eps \in (0,1/2)$ with $(1 +  \eps)^2 \tau^d
 \fmax <  \lambda_c$.
The Poisson  point process   $\Po_{n(1+\eps)}:=
\{X_1,\ldots,X_{M_{n(1+\eps)}}\}$, 
is stochastically dominated by $\H_{n \fmax (1+\eps)}$
(we say a point process $\X$ is stochastically dominated
by a point process $\Y$ if there exist coupled
point processes $\X', \Y'$ with $\X' \subset \Y'$
almost surely and $\X'$ having the distribution
of $\X$ and $\Y'$ having the distribution of $\Y$).
Hence by scaling, $\lambda^{1/d} \Po_{n(1+\eps)}$ is stochastically
dominated by $\H_{n \fmax (1+\eps)/ \lambda}$,
and hence we have for $n \leq \lambda (1+\eps) $ that
$\lla^{1/d} \Po_{n(1+\eps)}$ is stochastically
dominated by $\H_{ \fmax (1+\eps)^2 }$.  Therefore
for $u >0$, 
\bea
P[ R_{x,\lla,n} > u]
\leq 
P[ M_{n(1+\eps)} < n ] 
+ P[ R((1+ \eps)^2 \fmax) > u- \tau]. 
\lbl{0713a}
\eea
By scaling, the second probability 
in \eq{0713a} equals the probability that there is a path
from the origin in $\G(\H_{\tau^d (1+  \eps)^2 \fmax} \cup \{0\},1)$
to a point at Euclidean distance 
greater than $\tau^{-1}u -1$ from the origin.
By the exponential decay for subcritical continuum percolation,
(see e.g. Lemma 10.2 of \cite{Penbk}), 
this probability decays exponentially in $u$ (and does
not depend on $n$).

Let $\Delta:= \diam(\supp(f))$
(here assumed finite).
By Lemma \ref{PenL1.1}, the first term
in the right hand side of \eq{0713a} 
decays exponentially in $n$. Hence, 
there is a finite positive constant $C$, independent of $\lla$, such that
provided we have $n > (1-\eps) \lla^{1/d}$ we have
for all $ u \leq \lla^{1/d} (\Delta + \tau)$ that
$$
P[ M_{n(1+\eps)} < n ] \leq C \exp (- C^{-1} \lla^{1/d})  
\leq C \exp (- ((\Delta +\tau)C)^{-1} u) .
$$
On the other hand $P[ R_{x,\lla,n} > u]=0$
for $u > \lla^{1/d}(\Delta + \tau)$.  Combined with \eq{0713a} this shows that
there is a constant $C$ such that for all 
$(x,n,\lla,u) \in \supp(f) \times \N \times [1,,\infty)^2 $
with $n \leq (1+\eps) \lla$, we have
\bea
 P[ R_{x,\lla,n} > u]
\leq C \exp(-  u/C).
\lbl{0715a}
\eea
Now suppose $\A \in \TT_3 $, and $x \in \supp(f) $.
Let $R_{x,\lla,n,\A}$ be equal to $\tau$ plus
the maximum Euclidean  distance from $\lambda^{1/d}x$
of vertices in $\G(\lla^{1/d}(\X_n \cup \A \cup \{x\});\tau)$ 
that are pathwise connected to $\lambda^{1/d}x$.
Changes to $\X_n \cup \A$ at a Euclidean distance greater than
$\lambda^{-1/d} R_{x,\lambda,n,\A}$ 
from $x$  will have no effect on $\xi((\lambda^{1/d}x,T);
\lambda^{1/d}(\X_n^* \cup \A^*) )$; that is,
\eq{BiRS} holds.
To check the tail behaviour of $R_{x,\lla,n,\A}$,
suppose for example that $\A$ has three elements,
$x_1$, $x_2$  and $x_3$.
Then it is not hard to see that
$$
R_{x,\lla,n,\A} 
\leq 
R_{x,\lla,n} 
+ R_{x_1,\lla,n} 
+ R_{x_2,\lla,n} 
+ R_{x_3,\lla,n}, 
$$
and likewise  when $\A$ has fewer than three elements.
 Using this  together with \eq{0715a},  it is easy to
deduce that there is a constant $C$ such that 
for all
$(x,n,\A,\lla,u) \in \supp(f)  \times \N \times \TT_3
 \times [1,\infty)^2 $
with $n \leq (1+\eps) \lla$, and we have
\bea
 P[ R_{x,\lla,n,\A} > u]
\leq C \exp(-  u/C).
\lbl{0715a2}
\eea
In
other words, $\xi$ is  binomially exponentially stabilizing.

Next we check  the moments condition \eq{Bimoments},
with $p=4$ and using  the same choice of $\eps$ as before. 
By our definition of $\xi$ and the growth bound
\eq{polybd}, we have 
 for all $(x,n, \A, \lla) \in \supp(f) \times \N \times
  \TT_3 \times [1,\infty)^2 $  with $n \leq \lla(1+ \eps)$ 
that
\bea
\E[ \xi(( \lla^{1/d} x,T); \lla^{1/d} (\X_n^* \cup \A^*)) ^4]
\leq 
\polybeta^4 \E[ ( \card ( {\cal C}) + \diam ({\cal C}) )^{4\polybeta} ], 
\lbl{0715b}
\eea
where ${\cal C}$ is the vertex set of the
 component of $\G(\lla^{1/d}(\X_n \cup \A \cup \{x\});
\tau)$ containing $\lla^{1/d} x$.
By \eq{0715a2},  there is a constant $C$ such
that for all $(x,n, \A, \lla,u) \in (\supp(f) \times \N 
 \times \TT_3 \times [1,\infty)^2 $  with $n \leq \lla(1+ \eps)$ we have 
\bea
 P[ \diam ({\cal C}) >u ] \leq C \exp(-u/C);
\lbl{0715c}
\eea 
moreover,
\bea
P[
 \card ({\cal C}) >u ]  \leq 
 P[\diam ({\cal C}) >u^{1/(2d)} ] + 
P [ \card ( \X_n \cap B(x; \lla^{-1/d} u^{1/(2d)}))  > u-4] 
\nonumber \\
\label{0719a}
\eea
and the first term in  the right hand
side of \eq{0719a} decays exponentially in $u^{1/(2d)}$
by \eq{0715c}.
Since
$ \card ( \X_n \cap B(x; \lla^{-1/d} u^{1/(2d)})) $
is binomially distributed with
$$
 \E [ \card ( \X_n \cap B(x; \lla^{-1/d} u^{1/(2d)})) ] 
\leq u^{1/2} \omega_d \fmax n/\lla,
$$
by Lemma \ref{PenL1.1} there is a constant $C$ such that
 for all
$(x,n,\lla,u)$ with $n \leq \lla(1+\eps)$
we have that
$$
 P[ \card ( \X_n \cap B(x; \lla^{-1/d} u^{1/(2d)})) > u-4]
\leq C \exp (- C^{-1} u^{1/2} ).
$$ 
Thus by \eq{0719a} there is a constant, also denoted $C$, such that
 for all $(x,n,\A,\lla,u)$ with $n \leq \lla(1+ \eps)$ we have 
$$
 P[\card ({\cal C}) >u ]  \leq  C \exp ( -C^{-1} u^{-1/(2d)}),
$$ 
and combining this with \eq{0715c} and
using \eq{0715b} gives us a uniform tail bound 
which is enough to ensure \eq{Bimoments}.
The argument for \eq{Pomoments}  is similar.

Thus our $\xi$ satisfies all the assumptions of 
Theorem \ref{lemPEJP}, and we can deduce \eq{Hclteq} and
\eq{varconv} for some $\sigma \geq 0$ by applying 
that result with $\lla(n) = r_n^{-d}$.
  Then by applying Theorem   \ref{finranthm},
we can deduce that $\sigma >0$ and $\spann(H) < \infty$
 and \eq{1205a} holds whenever $\spann(H)|b$.
$\qed$ \\

{\bf Proof of  Theorem \ref{nnlclt}.}
Suppose the hypotheses of Theorem \ref{nnlclt}
hold, and
assume without loss of generality
that $\xi(\bx,\X^*)=0$ whenever $\X^* \setminus \{\bx\}$
has fewer than $\kappa$ elements.
We assert that under these hypotheses,
there exists a constant $C$ such that
 for all $(x,n,\lla,u) \in \supp(f) \times \N \times [1,\infty)^2$
with $n \in [\lla/2,3\lla/2]$
and  $n  \geq \kappa$, 
 we have
\bea
P[ \lla^{1/d} R_\kappa( ( x,T);
  \X_n^*) > u]  \leq C \exp (-C^{-1} u ).
\lbl{0714a}
\eea
Indeed, if $\supp(f)$ is a compact convex
region in $\R^d$ and $f$  is bounded away from zero on $\supp(f)$,
then \eq{0714a} is demonstrated in Section 6.3 of \cite{PenEJP},
while  if $ \supp(f)$ is
   a compact $d$-dimensional submanifold-with-boundary of
$\R^d$, and $f$ is bounded away from zero on $\supp(f)$,
then \eq{0714a} comes from the proof of
Lemma 6.1 of \cite{PYmfld}. 

It is easy to see that $\xi$ is
 $\lla$-homogeneously stabilizing for all $\lla >0$.
Also, for any $(x,\A) \in (\supp(f) \times \TT_3)$ we obviously have
$R_\kappa ((x,T); \X^* \cup \A^*) \leq
R_\kappa ((x,T); \X^* ) $ and hence by
\eq{0714a}, $\xi$ is  binomially exponentially 
stabilizing,
and exponential stabilization comes from a similar 
estimate with a Poisson sample.

We need  to check the moments conditions to be able
to deduce \eq{Hclteq} via
Theorem \ref{lemPEJP}.
With $\polybeta$ as in the growth bound \eq{polynbr}, we claim that
 there is a constant $C$ such that
for any ${\cal A} \in \TT_3$,
  any $x \in \supp(f)$,
 and any $u >0$, and
  for all $(x,n, \A ,\lla,u) \in  \supp(f) \times \N \times \TT_3
 \times [1,\infty)^2$ with
 $\lla/2 \leq n \leq 3\lla/2$,
and $n \geq \kappa$, we have
\bea
P [| \xi((\lla^{1/d} x,T);\lla^{1/d} (  
\X^*_n \cup {\cal A}^* ))| >  u]
\leq 
P [ \polybeta ( 1 + \lla^{1/d} R_\kappa((x,T),\X_n^*) )^\polybeta > u ]
\nonumber
\\
\leq C \exp( - C^{-1} u^{1/\polybeta}).
\lbl{0714b}
\eea
Indeed, the first bound comes from
the \eq{polynbr}, and the second bound comes from
\eq{0714a}. Using \eq{0714b}, 
 we can deduce the moments bound
\eq{Bimoments} for $p=4$ and $\eps =1/2$.
We can derive \eq{Pomoments} similarly.
Thus Theorem  \ref{lemPEJP}
is applicable, and enables us
to deduce \eq{Hclteq}  and \eq{varconv} for some $\sigma \geq 0$,
in the present setting.
Then by using Theorem \ref{finranthm},
 we can deduce that $\sigma >0$ and $\spann(H) >0$ and
\eq{1205a} holds whenever $\spann(H)|b$.
$\qed $  \\

 {\bf Acknowledgments.} We thank
 the Oberwolfach Mathematical Research Instiute
 for hosting
 the 2008 workshop `New Perspectives in Stochastic Geometry',
at which
this work was started. We also thank Antal J\'arai for helpful discussions.



\end{document}